\documentclass[letterpaper,11pt]{amsart}

\usepackage[margin=1.2in]{geometry}
\usepackage{amsmath,amsthm,amssymb, mathtools}
\usepackage{xspace,xcolor}
\usepackage[breaklinks,colorlinks,citecolor=teal,linkcolor=teal,urlcolor=teal,pagebackref,hyperindex]{hyperref}
\usepackage[alphabetic]{amsrefs}
\usepackage{comment}
\usepackage{enumerate}
\usepackage[all]{xy}
\usepackage{tikz-cd}
\usepackage{color}
\usepackage{comment}

\theoremstyle{plain}
\newtheorem{theo}{Theorem}[section]

\newtheorem{lemm}[theo]{Lemma}
\newtheorem{prop}[theo]{Proposition}

\theoremstyle{definition}
\newtheorem{defi}[theo]{Definition}

\newtheorem{eg}[theo]{Example}

\theoremstyle{remark}
\newtheorem{rema}[theo]{Remark}
\newtheorem{notation}[theo]{Notation}

\newcommand{\bfR}{\mathbf{R}}

\newcommand{\bb}[1]{\mathbb{#1}}

\newcommand{\QQ}{\mathbb{Q}}

\newcommand{\RR}{\mathbb{R}}
\newcommand{\CC}{\mathbb{C}}

\newcommand{\DD}{\mathbb{D}}

\newcommand{\ZZ}{\mathbb{Z}}

\newcommand{\cC}{\mathcal{C}}
\newcommand{\cD}{\mathcal{D}}

\newcommand{\cH}{\mathcal{H}}

\newcommand{\cM}{\mathcal{M}}
\newcommand{\cN}{\mathcal{N}}
\newcommand{\cO}{\mathcal{O}}

\newcommand{\uind}[1]{^{(#1)}}
\newcommand{\lind}[1]{_{(#1)}}
\newcommand{\ubind}[1]{^{[#1]}}

\newcommand{\lsta}{_{*}}

\newcommand{\sta}{^{*}}

\newcommand{\dual}{^{\vee}}
\newcommand{\ddual}{^{\vee\vee}}

\DeclareMathOperator{\Hom}{Hom}

\DeclareMathOperator{\Gr}{Gr}
\DeclareMathOperator{\gr}{gr}
\DeclareMathOperator{\an}{an}

\DeclareMathOperator{\Spec}{Spec}

\DeclareMathOperator{\SheafHom}{\mathcal{H}om}

\DeclareMathOperator{\coh}{coh}

\DeclareMathOperator{\pr}{pr}

\newcommand{\duBois}{\underline{\Omega}}
\newcommand{\DB}{\underline{\Omega}}

\newcommand{\de}{\partial}

\DeclareMathOperator{\DR}{DR}

\DeclareMathOperator{\IC}{IC}

\DeclareMathOperator{\MHM}{MHM}
\DeclareMathOperator{\HM}{HM}

\DeclareMathOperator{\dR}{dR}

\usepackage{soul}

\begin{document}
	\thanks{The first author was partially supported by NSF grant DMS-1952399.}
	
	\subjclass[2020]{14B05, 14C30, 14F10, 14M25, 14Q99,  32S35, 52B22}
	
	\author{Hyunsuk Kim}
	
	\address{Department of Mathematics, University of Michigan, 530 Church Street, Ann Arbor, MI 48109, USA}
	
	\email{kimhysuk@umich.edu}	

        \author{Sridhar Venkatesh}

        \address{Department of Mathematics, University of Michigan, 530 Church Street, Ann Arbor, MI 48109, USA}
	
	\email{srivenk@umich.edu}
	
	\begin{abstract} 
	       We study the Hodge filtration of the intersection cohomology Hodge module for toric varieties. More precisely, we study the cohomology sheaves of the graded de Rham complex of the intersection cohomology Hodge module and give a precise formula relating it with the stalks of the intersection cohomology as a constructible complex. The main idea is to use the Ishida complex in order to compute the higher direct images of the sheaf of reflexive differentials.
	\end{abstract} 
	
	\title[The intersection cohomology Hodge module of toric varieties]{The intersection cohomology Hodge module of toric varieties}
	
	\maketitle
	
	\setcounter{tocdepth}{1}
	\tableofcontents 
 
	
	\section{Introduction}

    A {\it toric variety} is a normal complex algebraic variety $X$ with an open subset isomorphic to the algebraic torus $(\mathbb{C}^*)^n$, along with an extension of the natural action of the torus to an action on $X$. Toric varieties provide an interesting interplay between algebraic geometry and convex geometry since they admit an alternate description in terms of convex geometric objects. As a consequence, algebro-geometric concepts on toric varieties correspond to much more elementary and tractable notions in convex geometry. One particular example where this relation is exploited is in studying intersection cohomology on toric varieties. There has been a long and fruitful study of the intersection cohomology complex and the intersection cohomology groups on toric varieties starting with the works of Stanley (\cite{Stanley:Intersection-cohomology-toric-varieties}) and Fieseler (\cite{Fieseler-ICprojtoric}), and more recently, the works of de Cataldo--Migliorini--Musta\c t\u a (\cite{dCMM-toricmaps}) and Saito (\cite{saito2020intersection}).

    However, the intersection cohomology complex has a richer structure as a (pure) Hodge module in the sense of Saito's theory (see \cite{saito1988modulesdeHodge}, \cite{saito1990mixedHodgemodules}). A (pure) Hodge module is a tuple $(\mathcal{M},F_\bullet,K,\alpha)$ satisfying certain conditions, where $\mathcal{M}$ is a holonomic $\cD$-module, $F_\bullet$ is a good filtration on $\cM$ (called the {\it Hodge filtration}), $K$ is a perverse sheaf on $X$ defined over $\mathbb{Q}$, and $\alpha$ is an isomorphism between $K \otimes_{\mathbb{Q}} \mathbb{C}$ and the analytic de Rham complex of $\cM$ 
    $$ \alpha \colon \DR^{\an}_{X}\cM \xrightarrow{\simeq} K \otimes_{\QQ} \CC.$$
    For more details, see Section \ref{section:hodge-modules}. Moreover, the Hodge filtration $F_{\bullet}$ on $\cM$ induces a natural filtration on $\DR_X(\mathcal{M})$, and the graded pieces $\gr_{k} \DR_{X}\cM$ lie inside the derived category of coherent sheaves on $X$.
    
    Given a variety $X$, it follows from Saito's theory that there exists a (pure) Hodge module $\IC_{X}^{H}$, whose underlying perverse sheaf is the intersection complex $\IC_X$. The main goal of this paper is to study the graded de Rham complex $\gr_k\DR_X \IC_X^{H}$\footnote{The shifted graded de Rham complex $\gr_{-p}\DR_X \IC_X^{H}[p-n]$ has been studied under the name \textit{intersection Du Bois complex}, denoted $I\DB^p_X$, in the recent papers \cites{park-popa, popa-shen-vo}.}, which can only be captured after enhancing $\IC_{X}$ to a Hodge module $\IC_{X}^{H}$. We now elaborate on how we study this object.
    
    Observe that any cohomology sheaf of the graded de Rham complex $\cH^{l}(\gr_{k} \DR_X \IC_{X}^{H})$ on an affine toric variety $X$ is an $\cO_{X}$-module, with a natural grading by $M$, where $M$ is the group of characters of the torus. We consider the generating function of $\dim_{\CC} (\cH^{l} \gr_{k} \DR_X \IC_{X}^{H})_{u}$ for $u \in M$. In order to compute this, we consider the following five generating functions for an $n$ dimensional affine toric variety $X$ defined by a cone $\sigma$, a proper birational toric morphism $\pi \colon Y \to X$ with $Y$ simplicial, and $\mu \subset \tau$ two faces of $\sigma$. These generating functions encode the following corresponding data:
    
    \begin{align*}
        \widetilde{F}_{\tau}(q) & = q^{-d_{\tau}} \sum_{j} h^{j}(\pi^{-1}(x_{\tau}), \QQ) q^{j} & \text{
        (Cohomology of fibers)} \\
        \widetilde{H}_{\mu, \tau}(q) & = q^{n-d_{\tau}}\sum_{j} h^{j}(\IC_{S_{\mu}})_{x_{\tau}} q^{j} &\text{(Intersection cohomology stalks)} \\
        D_{\tau}(q) & = \sum_{j} s_{\tau, j} q^{j} &\text{(Decomposition theorem)} \\
        \Omega_{\tau}(K,L) & = \sum_{k, l} \dim_{\CC} \left(R^{n+k+l}\pi\lsta \Omega_{Y}\ubind{-k}\right)_{u} K^{k} L^{l}, \quad u \in \tau_{\circ}\sta & \text{(Kähler differentials)} \\
        \dR_{\mu, \tau}(K, L) & = \sum_{k, l} \dim_{\CC} \left( \cH^{l} \gr_{k} \DR \IC^H_{S_{\mu}} \right)_{u} K^{k} L^{l}, \quad u \in \tau_{\circ}\sta & \text{(Graded de Rham complex).}
    \end{align*}
    For a detailed explanation of the notation, we refer to Section \ref{section:toric-var} for basic notation for toric varieties, Section \ref{section-Decompositiontheorem} for $\widetilde{F}_{\tau},\widetilde{H}_{\mu, \tau}$, and $D_{\tau}$, Section \ref{section:Kahler-diff} for $\Omega_{\tau}$, and Section \ref{section:grDR} for $\dR_{\mu, \tau}$. The relation between $\widetilde{F}_{\tau}$, $\widetilde{H}_{\mu, \tau}$, and $D_{\tau}$ is purely topological and well understood. The two extra pieces of data $\Omega_{\tau}$ and $\dR_{\mu, \tau}$ are related to the Hodge filtration on the intersection cohomology Hodge module $\IC_{X}^{H}$. However, we show that $\dR_{\mu, \tau}$ is completely determined by the combinatorial data of the toric variety by the following formula.

    \begin{theo}[Main Theorem = Theorem \ref{theorem:main-theorem-grDR}] \label{theo:main-theo1.1}
        With the above notation, we have
        $$ \dR_{\mu, \tau}(K, L) = \widetilde{H}_{\mu, \tau}(K^{-\frac{1}{2}}L) K^{\frac{d_{\mu}-d_{\tau}}{2}} (K^{-1} + L^{-1})^{n - d_{\tau}}. $$
    Moreover, $\widetilde H_{\mu, \tau}$ and hence $\dR_{\mu, \tau}$, can be computed explicitly in an algorithmic way.
    \end{theo}

    A rough sketch of the proof of the equality goes as follows. We know that $\widetilde{F}_{\tau}$, $\widetilde{H}_{\mu, \tau}$, and $D_{\tau}$ are related by the Decomposition theorem. Similarly, the decomposition theorem for Hodge modules gives a similar relation between the $\Omega_{\tau}$, $\dR_{\mu, \tau}$, and $D_{\tau}$. The main new ingredient is to compute $\Omega_{\tau}$ using the {\it Ishida complex}, which is presented in Section \ref{section:Kahler-diff}. Then we show that $\Omega_{\tau}$ can be expressed explicitly in terms of $\widetilde{F}_{\tau}$ when $\pi \colon Y \to X$ is given by a {\it barycentric subdivision} of the fan $\Sigma_{X}$ (Proposition \ref{prop-formulaKahlerdiff}). Along the way, we also show that barycentric subdivisions are shellable in Proposition \ref{prop:barycentric-subdivision-shellable}, which is an interesting combinatorial result in its own right. Finally, the rest of the proof follows using the two Decomposition theorems and induction.

    For the explicit calculation of $\widetilde H_{\mu, \tau}$ (and hence $\dR_{\mu, \tau}$), we again use the relation between $\widetilde{F}_{\tau}$, $\widetilde{H}_{\mu, \tau}$, and $D_{\tau}$ given by the Decomposition theorem. We have an explicit description of $\widetilde{F}_{\tau}$ from \cite{dCMM-toricmaps}. We then follow the strategy of \cite{CFS-Effectivedecompositiontheorem} to prove that once we have $\widetilde F_\tau$ explicitly, we can calculate \textbf{both} $\widetilde{H}_{\mu, \tau}$ and $D_{\tau}$ (Remark \ref{rema:explicit-computation-of-H}). We demonstrate this strategy in the appendix by explicitly computing $\widetilde{H}_{\mu,\tau}$ in dimensions $\leq 4$.

    \begin{rema}
        After the submission of our paper, we were informed by the referee of a possible alternative way to prove Theorem \ref{theo:main-theo1.1}. We properly spell this out in Appendix \ref{appendix-2}. \\
        \indent Although the approach in the appendix is a simpler way to prove Theorem \ref{theo:main-theo1.1}, we believe our approach yields several additional results along the way (e.g. developing a systematic way to calculate higher direct images of reflexive differential forms for toric morphisms, using shellability in these computations), which are interesting in their own right. In fact, we build on many of these ideas to provide a systematic understanding of the \textit{trivial Hodge module} $\QQ_{X}^{H}[n]$ on toric varieties, the local cohomology modules of toric varieties, and the Hodge structure on the singular cohomology of toric varieties in follow-up works \cite{LCDTV2, LCDTV1}.
    \end{rema}

    \section{Preliminaries}

    \subsection{Hodge modules} \label{section:hodge-modules}
    We give a brief summary on Hodge modules and state some results relevant to our situation. We will mostly follow the notation in \cite{saito1988modulesdeHodge} and \cite{saito1990mixedHodgemodules}. In these papers, Saito defines two abelian categories $\HM(X, w)$ and $\MHM(X)$ which are the categories of polarizable Hodge modules of weight $w$, and the category of polarizable mixed Hodge modules. The objects in $\HM(X, w)$ are holonomic $\cD$-modules with a filtration $F_{\bullet}$ by coherent sheaves with some extra structure satisfying suitable conditions. The objects in $\MHM(X)$ are also holonomic $\cD$-modules $\cM$ with a filtration $F_{\bullet}$ and an additional filtration $W_{\bullet}$ by holonomic $\cD$-modules satisfying some suitable conditions. The most important condition is that the graded piece $\gr_{w}^{W} \cM$ should be an object in $\HM(X, w)$. For general terminology associated to $\cD$-modules, we refer to \cite{HTT-Dmodulesbook}.
    
     The most important piece of data in our case is the filtration $F_{\bullet}$, known as the Hodge filtration. For a mixed Hodge module $\cM$ on a smooth variety $X$, we can consider the de Rham complex $$\DR_{X}\cM = [\cM \to \Omega_{X}^{1} \otimes \cM \to \ldots \to \Omega_{X}^{\dim X} \otimes \cM]$$
    which sits in cohomological degrees $-\dim X, \ldots, 0$. Moreover, $F_{\bullet}$ induces a natural filtration on $\DR_X(\mathcal{M})$ and the graded pieces $\gr_{k} \DR_{X}\cM$ are given by
    $$ \gr_{k} \DR_{X} \cM = [ \gr_{k} \cM \to \Omega_{X}^{1} \otimes \gr_{k+1} \cM \to \ldots \to \Omega_{X}^{\dim X} \otimes \gr_{k+ \dim X}\cM]$$
    which also sits in cohomological degrees $-\dim X, \ldots, 0$. The maps in this complex are $\cO_{X}$-linear and we thus view $\gr_{k} \DR_{X} \cM$ as an object in $\cD^b_{\rm coh}(X)$, the derived category of coherent sheaves of $X$. We mention that even if $X$ is singular, one can define the categories $\HM(X, w)$ and $\MHM(X)$ by embedding into a smooth variety $X \hookrightarrow Y$ and considering the objects in $\HM(Y, w)$ and $\MHM(Y)$, respectively, which are supported on $X$. The categories $\HM(X, w)$ and $\MHM(X)$ do not depend on the choice of the embedding. If $X$ cannot be embedded in a smooth variety, we locally embed each open set and impose suitable compatibility conditions on the intersections. The de Rham complex and the graded pieces $\gr_{k} \DR_{X}\cM$ are defined as above by locally embedding $X$ into a smooth variety. The graded pieces of the de Rham complex $\gr_{k} \DR_{X} \cM$ also do not depend on the choice of the embedding as objects in $\cD_{\coh}^{b}(X)$.

    We say a pure Hodge module $\cM \in \HM(X, w)$ has strict support $Z$ if $\cM$ is supported on $Z$ and has no nonzero subobjects or quotients supported on a strictly smaller subset of $Z$. The category $\HM(X, w)$ admits a decomposition by strict support, that means, for any $\cM \in \HM(X, w)$, there is a decomposition
    $$ \cM = \bigoplus_{Z \subset X} \cM_{Z}$$
    such that $\cM_{Z}$ has strict support $Z$, where the sum runs over all irreducible subvarieties of $X$. Also, a pure Hodge module $\cM \in \HM(X, w)$ with strict support $Z$ is a variation of Hodge structures $\cN$ on an open subset $U \subset Z$ of weight $w - \dim Z$. Conversely, any variation of Hodge structures on an open subset of the smooth locus of $Z$ can be uniquely extended to a pure Hodge module on $X$ with strict support $Z$. In this case, the underlying $\cD$-module is the intermediate extension of the $\cD$-module corresponding to the variation of Hodge structures.
    In this sense the intersection cohomology $\cD$-module on $X$ underlies a pure Hodge module of weight $\dim X$ since it is associated to the trivial variation of Hodge structures on the smooth locus of $X$. We denote by $\IC_{X}^{H}$ the intersection cohomology Hodge module of $X$ in order to distinguish this from the perverse sheaf $\IC_{X}$. 

    The derived category of mixed Hodge modules $\cD^{b} \MHM(X)$ has a six functor formalism and moreover, these functors are compatible with the functors at the level of perverse sheaves. The most important functor that we use is the pushforward $ \pi\lsta$, where $\pi$ is a proper morphism. We recall Saito's decomposition theorem for Hodge modules. We say a complex $\cM \in \cD^{b}(\MHM(X))$ is pure of weight $w$ if $\cH^{j} \cM$ is a pure Hodge module of weight $w + j$ for all $j$. In this case, we always have the following
    \begin{prop}[\cite{saito1990mixedHodgemodules}*{4.5.4}]
        If $\cM \in \cD^{b}(\MHM(X))$ is pure of weight $w$, then
        $$ \cM \simeq \bigoplus_{j} \cH^{j}(\cM) [-j].$$
    \end{prop}

    Similarly, we say that $\cM \in \cD^{b} \MHM(X)$ is of weight $\leq n$ (resp. $\geq n$), if the following condition is satisfied:
    $$ \Gr_{i}^{W} \cH^{j} \cM = 0 \quad \text{for} \quad i > n + j \quad (\text{resp. } i < n + j).$$
    By \cite{saito1990mixedHodgemodules}*{4.5.2}, if $\cM$ is of weight $\leq n$ (resp. $\geq n$), then $f_{!}\cM$ (resp. $f\lsta \cM$) is also of weight $\leq n$ (resp. $\geq n$). Since $f_{!} = f\lsta$ for proper morphisms, in this case $f\lsta$ takes pure complexes to pure complexes. Therefore, we have the following decomposition theorem.

    \begin{theo}[Saito's decomposition theorem] \label{theo-Saitodirectimage}
        Let $\pi \colon Y \to X$ be a proper morphism and $\cM \in \HM(Y, w)$ be a polarizable pure Hodge module. Then we have a decomposition
        $$ \pi\lsta \cM \simeq \bigoplus_{j \in \ZZ} \cH^{j} \pi\lsta \cM [-j] $$
        and $\cH^{j}\pi\lsta \cM \in \HM(X, w + j)$ for all $j$.
    \end{theo}

    We also recall that taking the graded de Rham complex commutes with the pushforward by a proper morphism $\pi \colon Y \to X$ \cite{saito1988modulesdeHodge}*{2.3.7}:
    \begin{equation} \label{equation-grDRcommuteswithpropermap}
       \gr_{k} \DR_{X}  \pi\lsta \cM \simeq \bfR \pi\lsta (\gr_{k} \DR_{Y} \cM) . 
    \end{equation}

    We end this section by discussing the relation between the Du Bois complex and mixed Hodge modules. In \cite{DuBois:complexe-de-deRham}, Du Bois introduced a filtered complex $\duBois_{X}^{\bullet}$ which can be thought of as a replacement of the de Rham complex $\Omega_{X}^{\bullet}$ when $X$ is singular. By taking the graded quotients, the {\it $p$-th Du Bois complex} is defined as
    $$ \duBois_{X}^{p} := \gr_{F}^{p} \duBois_{X}^{\bullet} [p].$$
    We have a natural comparison map $\Omega_{X}^{p} \to \duBois_{X}^{p}$ which is an isomorphism if $X$ is smooth. Note that $\duBois_{X}^{p}$ is an object in $\cD^b_{\rm coh}(X)$.

    In \cite{Saito-MixedHodgecomplexes}, Saito gives a description of the Du Bois complex using the {\it trivial mixed Hodge module} $\QQ_{X}^{H}$. The category of mixed Hodge modules over a point can be identified with the category of mixed Hodge structures. Hence, we have the Hodge module $\QQ_{\mathrm{pt}}^{H}$ with weight zero given by the following mixed Hodge structure $(V, W_{\bullet}, F^{\bullet})$:
    $$V = \QQ, \quad W_{-1}V = 0, \quad W_{0}V = V,\quad  F^{0}V_{\CC} = V_{\CC}, \quad F^{1}V_{\CC} = 0. $$
    For an arbitrary variety $X$, $\QQ_{X}^{H}$ is defined as
    $$ \QQ_{X}^{H} := (a_{X})\sta \QQ_{\mathrm{pt}}^{H} \in \cD^{b} (\MHM(X)),$$
    where $a_{X} \colon X\to \{ \mathrm{pt}\}$ is the structure morphism. It is a consequence of \cite{Saito-MixedHodgecomplexes}*{Theorem 4.2} that the graded de Rham complex of this Hodge modules is related to the Du Bois complex in the following way:
    $$ \gr_{k}^{F} \DR_{X} \QQ_{X}^{H}[n] \simeq \duBois_{X}^{-k} [n+k].$$
    One can also get this easily using \cite{Mustata-Popa:localcohomologyHodge}*{Proposition 5.5} and duality.

    The two objects $\QQ_{X}^{H}$ and $\duBois_{X}^{p}$ have nice descriptions when the variety $X$ has quotient singularities. In general, we have a natural morphism $\QQ_{X}^{H}[n] \to \IC_{X}^{H}$ in the derived category of mixed Hodge modules (see \cite[4.5.11]{saito1990mixedHodgemodules}). If $X$ has quotient singularities, this is an isomorphism at the level of perverse sheaves by \cite{borho-macpherson}*{Section 1.4}, which implies that $\QQ_{X}^{H}[n] \to \IC_{X}^{H}$ is also an isomorphism of mixed Hodge modules. Moreover, $\duBois_{X}^{p}$ coincides with the reflexive Kähler differentials $\Omega_{X}\ubind{p} := (\Omega_{X}^{p})\ddual$ if $X$ has quotient singularities \cite{DuBois:complexe-de-deRham}*{Théorème 5.3}. Wrapping all up, we have the following lemma:

    \begin{lemm} \label{lemm:quotientsing-grDR}
        If $X$ has quotient singularities, then
        $$ \gr_{k}^{F} \DR_{X} \IC_{X}^{H} \simeq \Omega_{X}\ubind{-k}[n+k]. $$
    \end{lemm}

    \subsection{Toric varieties} \label{section:toric-var}  

    Fix a free abelian group $N$ of rank $n$ and let $M := \Hom_{\mathbb{Z}}(N,\mathbb{Z})$. Denote $N_{\mathbb{R}} := N \otimes \mathbb{R}$ and $M_{\RR}:= M \otimes \mathbb{R}$. To a strongly convex rational polyhedral cone $\sigma \subset N_{\mathbb{R}}$, we associate an $n$-dimensional affine toric variety $X_\sigma$. More generally, to a fan $\Delta \subset N_{\mathbb{R}}$, we associate an $n$-dimensional toric variety $X_\Delta$ by gluing the affine toric varieties corresponding to the cones of $\Delta$. For general notions regarding toric varieties, we refer to \cite{Fulton-ToricVar} and \cite{CoxLittleSchenck-ToricVar}.

    \noindent {\bf Notation and terminology.} We collect some notation for convex cones that we will use. Here, $\sigma \subset N_{\mathbb{R}}$ denotes a strongly convex rational polyhedral cone, and $\mu, \tau$ denote faces of $\sigma$.
    \begin{enumerate}
        \item $\tau^\vee := \{ u \in M_{\mathbb{R}} \mid u(v) \geq 0, \forall v \in \tau \}$
        \item $\tau^{\perp} := \{ u \in M_{\RR} \mid u(v) = 0, \forall v \in \tau\}$
        \item $\tau\sta := \tau^{\perp} \cap \sigma\dual$
        \item $\tau_{\circ}\sta := (\tau\sta \cap M) \setminus \left( \bigcup_{\tau\subsetneq \nu} \nu\sta \right)$.
        \item $\langle \tau \rangle \subset N_{\RR}$ is the subspace spanned by $\tau$.
        \item $d_\tau := \dim_{\mathbb{R}} \langle \tau \rangle$.
        \item We denote by $\bar{\tau}_{\mu}$ the image of $\tau$ under the projection $N_{\RR} \to N_{\RR} / \langle \mu \rangle$.
        \item $\sigma$ is {\it full-dimensional} if $\langle \sigma \rangle = N_{\RR}$.
        \item $\sigma$ is {\it simplicial} if the 1-dimensional faces (i.e. rays) are linearly independent over $\RR$ in $N_{\RR}$.
    \end{enumerate}
    
    We collect several facts on toric varieties that we will need.
    \begin{rema} \label{rema:Basics-Toric-varieties}
    \begin{enumerate}
        \item \cite{CoxLittleSchenck-ToricVar}*{Theorem 9.2.5} For a proper toric morphism $\pi \colon Y\to X$, we have $R^{p}\pi\lsta \cO_{Y} = 0$ for all $p > 0$.

        \item \cite{Fulton-ToricVar}*{Section 3.1} Let $X_{\sigma} = \Spec \CC [\sigma\dual \cap M]$ be the affine toric variety corresponding to a strongly convex rational polyhedral cone $\sigma$. For an $r$-dimensional face $\tau \subset \sigma$, we get a torus invariant subvariety $S_{\tau} = \Spec \CC[\sigma\dual \cap \tau^{\perp} \cap M] \subset X_{\sigma}$ of codimension $r$. This is the affine toric variety corresponding to the cone $\bar{\sigma}_{\tau}$, where the lattice and the dual lattice are given by
        $$ N_{\tau} := \frac{N}{N \cap \langle\tau\rangle},\qquad M_{\tau}: = M \cap \tau^{\perp}.$$
        We denote by $O_{\tau} = \Spec \CC[M_{\tau}]$ the torus orbit corresponding to $\tau$. Also, $U_{\tau} = \Spec\CC[\tau\dual \cap M]$ is an open subset of $X_{\sigma}$ and we have the diagram of torus equivariant morphisms
        $$ \begin{tikzcd}
            U_{\tau} \ar[r, hook] \ar[d] & X_{\sigma} \ar[d] \\ O_{\tau} \ar[r, hook]  & S_{\tau}.
        \end{tikzcd} $$
        After fixing a non-canonical splitting $N = N_{\tau} \oplus (N \cap \langle\tau\rangle)$ and the corresponding splitting $M = M_{\tau} \oplus M'$, we can identify $U_{\tau} \to O_{\tau}$ as the projection $U_{\tau} = V_{\tau} \times O_{\tau} \to O_{\tau}$, where $V_{\tau}$ is the full-dimensional toric variety $\Spec \CC[\tau\dual \cap M']$.
        
        For two faces $\mu \subset \tau$, we denote by $V_{\mu, \tau}$ the full-dimensional affine toric variety corresponding to the cone $\bar{\tau}_{\mu} \subset \langle \tau \rangle/ \langle \mu \rangle$. We have an analogous diagram
        $$ \begin{tikzcd}
            S_{\mu} \ar[d] & U_{\bar{\tau}_\mu} \simeq V_{\mu, \tau} \times O_{\tau} \ar[d] \ar[l, hook']  \\
            S_{\tau} & \ar[l, hook'] O_{\tau}.
        \end{tikzcd} $$
        Note that
        \begin{align*}
            \Sigma_{S_{\mu}} & = \{ \bar{\tau}_{\mu} \colon\tau \in \Sigma_{X},  \mu \subset \tau \},\\
            \Sigma_{V_{\tau}} & = \{ \mu \in \Sigma_{X} \colon \mu \subset \tau \}, \\
            \Sigma_{V_{\mu, \tau}} &= \{ \bar{\nu}_{\mu} \colon \nu \in \Sigma_{X} , \mu \subset \nu \subset \tau \}
        \end{align*}
        as a fan in $N/\langle \mu \rangle \cap N$, $\langle\tau\rangle \cap N$, and $\langle \tau \rangle \cap N/ \langle \mu \rangle\cap N$, respectively.
        
        \item We say that a toric variety $X$ is simplicial if all the cones $\sigma$ in the fan are simplicial. If $X$ is simplicial, then $X$ has quotient singularities \cite[Theorem 11.4.8]{CoxLittleSchenck-ToricVar}. By Lemma \ref{lemm:quotientsing-grDR}, in this case we have a canonical isomorphism
        $$ \gr_{k}\DR_{X} \IC_{X}^{H} \simeq \Omega_{X}\ubind{-k}[n+k].$$

        \item \cite[Lemma 3.5]{oda-ConvexBodies} (see also \cite{Ishida}) There is a natural resolution of the sheaf of reflexive differentials $\Omega_{Y}\ubind{p}:= (\Omega_{Y}^{p})\ddual$ of a simplicial toric variety $Y$, called the \textit{Ishida complex}, that we now recall. Let $\tau \in \Sigma_{Y}$ be an $r$-dimensional cone. Since $Y$ is simplicial, $\tau$ is generated by $r$ rays $\rho_{1},\ldots, \rho_{r}$. We set
        $$ V_{\tau}^{p} := \left( \bigwedge^{p-r} \tau^{\perp} \right) \otimes \frac{M_{\CC}}{\rho_{1}^{\perp}} \otimes \ldots \otimes \frac{M_{\CC}}{\rho_{r}^{\perp}}, $$
        where $M_{\CC} := M \otimes_{\ZZ} \CC$. For instance, $V_{0}^{p} = \bigwedge^{p} M_{\CC}$. Then we have an exact complex
        \begin{equation}\label{equation:Ishida-complex}
            0 \to \Omega_{Y}\ubind{p} \to V_{0}^{p} \otimes \cO_{Y} \to \bigoplus_{\substack{\tau \in \Sigma_{Y} \\ \dim \tau = 1}} V_{\tau}^{p} \otimes \cO_{S_{\tau}} \to \bigoplus_{\substack{\tau \in \Sigma_{Y} \\ \dim \tau = 2}} V_{\tau}^{p} \otimes \cO_{S_{\tau}} \to \ldots \to \bigoplus_{\substack{\tau \in \Sigma_{Y} \\ \dim \tau = p}} V_{\tau}^{p} \otimes \cO_{S_{\tau}}\to 0.
        \end{equation}
    \end{enumerate}
    \end{rema}
    
    The relevant situation for us will be when $\pi \colon Y \to X$ is obtained by a simplicial subdivision of the fan $\Sigma_{X}$. In this case, for $a \in \Sigma_{Y}$, we denote by $\pi\lsta (a)$ the minimal cone of $\Sigma_{X}$ containing $a$. For $\tau \in\Sigma_{X}$, we let
    $$ d_{l}(\tau) = \# \{ a \in \Sigma_{Y} \colon \dim(a) =l, \pi\lsta(a) = \tau \}.$$
    We point out that our notation for $d_{l}(\tau)$ is slightly different from the one in \cite{dCMM-toricmaps}. The following proposition describes how the fibers look like for arbitrary proper toric morphisms.

    \begin{prop}[\cite{dCMM-toricmaps}*{Lemma 2.6, Proposition 2.7}] \label{prop-fibertoric}
        Let $\pi \colon Y \to X$ be a proper toric morphism.
        \begin{enumerate}[(i)]
            \item Then every irreducible component of the fiber $\pi^{-1}(x)$ is a toric variety. Moreover, this is smooth (resp. simplicial) if $Y$ is smooth (resp. simplicial).
            \item For any $x_\tau \in O_\tau$, we have an isomorphism $\pi^{-1}(O_\tau) \simeq \pi^{-1}(x_\tau) \times O_\tau$ such that the restriction of $\pi$ to $\pi^{-1}(O_\tau)$ corresponds to the projection onto the second component. In particular, $\pi^{-1}(x'_\tau) \simeq \pi^{-1}(x_\tau)$ for every $x'_\tau \in O_\tau$.
        \end{enumerate}
    \end{prop}
    
    The following proposition gives a combinatorial formula for the cohomology of the fibers.

    \begin{prop}[\cite{dCMM-toricmaps}*{Theorem B, Corollary C}] \label{prop-fibercohomology}
        Let $\pi \colon Y \to X$ be the proper birational toric morphism obtained by a simplicial subdivision of $\Sigma_{X}$. For every $\tau \in \Sigma_{X}$ and every $x_{\tau} \in \cO_{\tau}$, the Hodge structure on $H^{j}(f^{-1}(x_{\tau}), \QQ)$ is pure of Hodge--Tate type. Moreover, we have the following formula
        $$ \sum_{j} \dim_{\QQ} H^{j}(f^{-1}(x_{\tau}), \QQ) \cdot q^{j} = \sum_{l} d_{l}(\tau) \cdot (q^{2}-1)^{d_{\tau} - l}.$$
    \end{prop}

    \begin{rema} \label{rema:barycentric-resolution}
        At last, we describe the construction of what we call a {\it barycentric resolution} $Y$ of an affine toric variety $X$ of dimension $n$. Let $X$ be the affine toric variety corresponding to a cone $\sigma$. Since $\sigma$ is strictly convex, there exists a linear functional $l \in M_{\mathbb{Q}}$ such that $l \geq 0$ on $\sigma$ and $\{l = 0\}\cap \sigma = \{0\}$. Then consider the polytope $P = \sigma \cap \{ l = 1\}$. Note that the cone generated by $P$ is $\sigma$. For each face $\mu \subset \sigma$, choose $\rho_{\mu} \in N_{\QQ} \cap P$ lying inside the relative interior of $\mu$. We construct the fan $\Sigma_{Y}$ by describing its maximal cones. Each maximal cone in $\Sigma_{Y}$ is of the form
        $$ \tau_{\mu_{1},\ldots, \mu_{n}} = \mathrm{span}_{\RR_{\geq 0}} \{ \rho_{\mu_{1}}, \ldots, \rho_{\mu_{n}} \} ,$$
        where each $\mu_{i}$ is an $i$-dimensional face of $\sigma$ satisfying $\mu_{1}\subset \mu_{2} \subset \ldots \mu_{n} = \sigma$. Note that $\Sigma_{Y}$ is simplicial by construction. We have a proper birational toric morphism $\pi \colon Y \to X$ corresponding to the map $\Sigma_{Y} \to \Sigma_{X}$ induced by the identity on $N$. We call $\Sigma_Y$ a {\it barycentric subdivision} of $\sigma$ and the corresponding toric morphism $\pi\colon Y \to X$ a {\it barycentric resolution} of $X$. Observe that we also get a simplicial polytopal complex $\cC_Y$ (see Definition \ref{defi:polytopal-complex}), the cone over which is $\Sigma_Y$. Each maximal simplex of $\cC_Y$ corresponds to a sequence of faces
        $$ P = F\uind{n}, F\uind{n-1}, F\uind{n-2}, \ldots , F\uind{0}$$
        where $F\uind{i}$ is a facet of $F\uind{i+1}$ for $0 \leq i \leq n-1$. Observe that $F\uind{0}$ is a vertex of $P$. Given this sequence, the set of vertices of the corresponding simplex consists of
        $$\rho_{\mathrm{span}_{\mathbb{R}_{\geq 0}}F\uind{n}}, \rho_{\mathrm{span}_{\mathbb{R}_{\geq 0}}F\uind{n-1}},\ldots, \rho_{\mathrm{span}_{\mathbb{R}_{\geq 0}}F\uind{1}} ,\rho_{\mathrm{span}_{\mathbb{R}_{\geq 0}}F\uind{0}}.$$    
        Even though the toric variety $Y$ constructed in this way depends on the choice of the generators in the relative interior, this will not affect the arguments throughout this article.
    \end{rema}

    \begin{rema}
        Though it will not be important for what follows, we mention that every barycentric resolution $\pi \colon Y \to X$ of an affine toric variety is a projective morphism.
    \end{rema}

    \subsection{Polytopes and Shellability}\label{section-shellability}
    In this section we review the concept of shellability, as we will later use it to prove results about the pushforward of the Ishida complex. We follow \cite[Chapter 8]{ziegler-polytopes}, where the reader can find additional details.

    \begin{defi}\cite{ziegler-polytopes} \label{defi:polytopal-complex}
        A {\it polytopal complex} is a finite, non-empty collection $\cC$ of polytopes (called faces of $\cC$) in $\RR^N$ that contains the faces of all its polytopes and such that the intersection of any two of its polytopes is a face of each of them. The inclusion-maximal faces of $\cC$ are called the \textbf{facets} of $\cC$.

        A polytopal complex $\cC$ is {\it pure} if all its facets have the same dimension and is {\it simplicial} if all its faces are simplices.
    \end{defi}

    \begin{eg}
        If $P$ is a polytope, then the boundary complex $\cC(\partial P)$, which is defined to be the set of all proper faces of $P$, is a pure polytopal complex of dimension $\dim(P)-1$.
    \end{eg}

    \begin{defi}[\cite{ziegler-polytopes}*{Definition 8.1}]
        Let $\cC$ be a pure $d$-dimensional polytopal complex. A {\it shelling} of $\cC$ is a linear ordering $F_1,F_2,\dots,F_s$ of the facets of $\cC$ such that either $\cC$ is a set of points, or it satisfies the following condition:
        \begin{enumerate}
            \item The boundary complex $\cC(\de F_{1})$ of the first facet $F_{1}$ has a shelling.

            \item For $1 < j \leq s$,
            $$ F_{j} \cap \left( \bigcup_{i=1}^{j-1} F_{i} \right) = G_{1} \cup G_{2} \cup \ldots \cup G_{r}$$
            for some shelling $G_{1},G_{2}\ldots, G_{r},\ldots, G_{t}$ of $\cC(\de F_{j})$
        \end{enumerate}

        A pure polytopal complex is {\it shellable} if it has a shelling.
    \end{defi}

    We will use the following theorem of Bruggesser and Mani:

    \begin{theo}[\cite{Bruggesser-Mani:Shellable}]\label{thm:bruggesser-mani}
        For a polytope $P$, the polytopal complex $\cC(\de P)$ is shellable.
    \end{theo}

    \begin{defi}[Type] \label{defi:type-of-shelling}
        Let $F_1,\dots,F_r$ be a shelling of a simplicial polytopal complex. Define $F_1$ to be of type $0$. For $j \geq 2$, we have that $F_j \cap (\bigcup_{i=1}^{j-1} F_i)$ is a pure $(d-1)$-dimensional complex. Define $F_j$ to be of {\it type} $l$, where $l$ is the number of facets in the pure complex $F_j \cap (\bigcup_{i=1}^{j-1} F_i)$.
    \end{defi}

    \begin{notation}
        Let $F$ be a simplex whose vertices are $v_{1},\ldots, v_{n}$. In this case, we sometimes use the notation 
        $$ F = [v_{1},\ldots, v_{n}]$$
        in order to denote $F$.
    \end{notation}
    
    \begin{rema} \label{rema:type-l-when-add}
        With notation as in Definition \ref{defi:type-of-shelling}, let $F_{j} = [v_{1},\ldots, v_{n}]$ be of type $l$, and let $G_{k} = [v_{1}, \ldots , \hat{v}_{j_{k}}, \ldots, v_{n}]$ for $k= 1,\dots,l$ be the facets in the pure complex $F_j \cap (\bigcup_{i=1}^{j-1} F_i)$, so that 
        \[ F_j \cap \left(\bigcup_{i=1}^{j-1} F_i\right) = \bigcup_{k=1}^{l} G_k.\]
        Observe that we have:
        \begin{align*}
            \left\{\text{Faces of $F_j$ containing $[v_{j_{1}},\ldots, v_{j_{l}}]$}\right\} &= \left\{ \text{Faces of $F_j$ not contained in any $G_k$ $\forall k=1,\dots,l$} \right\}\\
            &= \left\{ \text{Faces of $F_j$ not contained in $\bigcup_{k=1}^lG_k$} \right\}\\
            &= \left\{ \text{Faces of $F_j$ not contained in $F_j \cap \left(\bigcup_{i=1}^{j-1} F_i\right)$} \right\}.
        \end{align*}
    \end{rema}

    We will be interested in the shellability of the pure simplicial polytopal complex $\cC_Y$ as defined in Remark \ref{rema:barycentric-resolution}.
    
    \begin{prop}\label{prop:barycentric-subdivision-shellable}
        With notation as in Remark \ref{rema:barycentric-resolution}, the polytopal complex $\cC_Y$ is shellable.
    \end{prop}

    \begin{proof}
        We describe a shelling of $\cC_Y$ by defining a lexicographic order on the set of maximal simplices of $\cC_Y$.
        
        Let us recall the maximal simplices of the simplicial polytopal complex $\cC_{Y}$. Each maximal simplex corresponds to a sequence of faces
        $$ P = F\uind{n}, F\uind{n-1}, F\uind{n-2}, \ldots , F\uind{0}$$
        where $F\uind{i}$ is a facet of $F\uind{i+1}$ for $0 \leq i \leq n-1$. Observe that $F\uind{0}$ is a vertex of $P$. From now on, for a maximal simplex $\Delta$ in $\cC_{Y}$ corresponding to the chain of faces $P = F\uind{n},\ldots, F\uind{0}$, we use the notation
        $$ \Delta = (F\uind{n} \supset \ldots \supset F\uind{0}) $$
        to represent the chain of faces.

        First, let $F_1,\dots,F_r$ be an ordering of the facets of the polytope $P$ such that it gives a shelling of $\cC(\partial P)$ (such an ordering exists by Theorem \ref{thm:bruggesser-mani}). This defines an ordering, call it $\prec_{F\uind{n}}$, on the set of facets of $P =F\uind{n}$, given by:
        \[ F_a \prec_{F\uind{n}} F_b \text{ for } a<b.  \]
        Similarly, given an $n-1$ dimensional face $F\uind{n-1}$ of $P$, let us now define an order $\prec_{F\uind{n},F\uind{n-1}}$ on the set of facets of $F\uind{n-1}$. Observe that
        \[F\uind{n-1} \cap \left(\bigcup_{G\uind{n-1} \prec_{F\uind{n}} F\uind{n-1}} G\uind{n-1} \right) = F\uind{n-2}_1 \cup \dots \cup F\uind{n-2}_l , \]
        where $F\uind{n-2}_1 , \ldots , F\uind{n-2}_l , \ldots , F\uind{n-2}_r$ is a shelling of $\partial F\uind{n-1}$. Define:
        \[ F\uind{n-2}_a \prec_{F\uind{n},F\uind{n-1}} F\uind{n-2}_b \text{ for $a<b$}.\]
        In the same vein, given a chain $F\uind{n} \supset \dots \supset F\uind{i}$, we now inductively define an order $\prec_{F\uind{n},\dots,F\uind{i}}$ on the set of facets of $F\uind{i}$. Observe that, by induction, the chain $F\uind{n} \supset \dots \supset F\uind{i+1}$ defines an order on the set of facets of $F\uind{i+1}$ (by defining a shelling on $F\uind{i+1}$). We will use this information to get an order $\prec_{F\uind{n},\dots,F\uind{i}}$ on the set of facets of $F\uind{i}$. Observe:
        \[ F\uind{i} \cap \left(\bigcup_{G\uind{i} \prec_{F\uind{n},\dots,F\uind{i+1}} F\uind{i}} G\uind{i} \right) = F\uind{i-1}_1 \cup \dots \cup F\uind{i-1}_l, \]
        where $F\uind{i-1}_1 , \ldots , F\uind{i-1}_l , \ldots , F\uind{i-1}_r$ is a shelling of $\partial F\uind{i}$. Define:
        \[ F\uind{i}_a \prec_{F\uind{n},\dots,F\uind{i}} F\uind{i}_b \text{ for $a<b$}.\]

        We can now finally define the lexicographic order on the set of all simplices of $\cC_Y$: Given two simplices $\Delta = (F\uind{n} \supset \ldots \supset F\uind{0})$ and $\Delta' = (G\uind{n}\supset \ldots \supset G\uind{0})$, let $k =\max\{ i \mid F\uind{i} \neq G\uind{i}\}$ (observe that $k<n$ as $F\uind{n} = G\uind{n}$), then define $\Delta \prec \Delta'$ if
        \[ F\uind{k} \prec_{F\uind{n},\dots,F\uind{k+1}} G\uind{k}. \]
        Now, we finally prove that if we arrange the simplices in the lexicographical order, this is a shelling of $\cC_Y$. Consider a simplex
        $$\Delta = (F\uind{n}\supset \ldots \supset F\uind{0})$$
        and write the vertices of $\Delta$ as $v\lind{n}, v\lind{n-1}, \ldots, v\lind{0}$ where $v\lind{i} = \rho_{\mathrm{span}_{\mathbb{R}_{\geq 0}}F\uind{i}}$, with notation as in Remark \ref{rema:barycentric-resolution}. One observation is the following: if we have a chain of faces $F\uind{i+2} \supset F\uind{i+1} \supset F\uind{i}$, then there is a unique face $F_{\ast}\uind{i+1} \neq F\uind{i+1}$ such that $F\uind{i+2} \supset F\lsta\uind{i+1} \supset F\uind{i}$. Denote $F\lsta\uind{0}$ as the unique vertex of $F\uind{1}$ which is not $F\uind{0}$.
        
        We describe the simplices which share a facet with $\Delta$. There are exactly $n$ of them:
        \begin{align*}
           \Delta_{i} &= (F\uind{n}\supset \ldots \supset F\uind{i+1} \supset F\uind{i}\lsta \supset F\uind{i-1} \supset \ldots \supset F\uind{0}), \qquad \text{for } 0\leq i \leq n-1.
        \end{align*}
        Let $i_{1}> \ldots > i_{l}$ be the integers such that $\Delta_{i_{t}} \prec \Delta$ exactly for these indices.

        Note that the facets $[v\lind{n},\ldots, \widehat{v\lind{i_{t}}}, \ldots, v\lind{0}]$ are the facets of $\Delta$ which are contained in $\bigcup_{\Delta' \prec \Delta}\Delta'$. We claim that any face of $\Delta$ containing $[v\lind{i_{1}},\ldots, v\lind{i_{t}}]$ is not a face of $\Delta'$ with $\Delta' \prec \Delta$. Provided that, it is straightforward to see that
        $$ \Delta \cap \bigcup_{\Delta' \prec \Delta} \Delta' = \bigcup_{t=1}^{l} [v\lind{n},\ldots, \widehat{v\lind{i_{t}}}, \ldots, v\lind{0}],$$
        following Remark \ref{rema:type-l-when-add}. This also shows that $\Delta$ is of type $l$ and that $\prec$ is a shelling order.
        
        We are left with proving that any face of $\Delta$ containing $[v\lind{i_{1}},\ldots, v\lind{i_{t}}]$ is not a face of $\Delta'$ with $\Delta' \prec \Delta$, that is, we want to prove that the set $S$ of simplices $\Delta' =( G\uind{n} \supset \ldots \supset G\uind{0} )$
        satisfying the following properties:
        \begin{enumerate}
            \item $\Delta' \prec \Delta$
            \item $G\uind{i_{t}} = F\uind{i_{t}}$ for $t = 1, \ldots, l$
        \end{enumerate}
        is empty.
        
        Assume for the sake of contradiction that $S \neq \emptyset$. For each $\Delta' =( G\uind{n} \supset \ldots \supset G\uind{0} )  \in S$ consider
        $$ i_{\Delta'} = \max \{ i : G\uind{i} \neq F\uind{i} \}.$$
        Pick $\Delta'$ such that $i = i_{\Delta'}$ is the minimum among these quantities. We note that $i$ cannot be one of the $i_{1},\ldots, i_{l}$ by the definition of $S$.
        
        First, we suppose that $i > i_{l}$. We have $G\uind{i} \prec_{F\uind{n},\ldots, F\uind{i+1}} F\uind{i}$. Also, we have $F\uind{i_{l}} = G\uind{i_{l}}$ which implies that $F\uind{i} \cap G\uind{i}$ contains $F\uind{i_{l}}$. This implies that there exists a facet $\widetilde{G}\uind{i} \prec_{F\uind{n},\ldots, F\uind{i+1}} F\uind{i}$ of $F\uind{i+1}$ such that $\widetilde{G}\uind{i-1} :=  \widetilde{G}\uind{i} \cap F\uind{i}$ is a $(i-1)$-dimensional face and $F\uind{i_{l}}$ is a face of $\widetilde{G}\uind{i-1}$. We replace $\Delta'$ with
        $$ \Delta'' = (G\uind{n} \supset \ldots \supset F\uind{i+1} \supset \widetilde{G}\uind{i} \supset \widetilde{G}\uind{i-1}\supset \ldots \supset F\uind{i_{l}} \supset \ldots \supset F\uind{0}). $$
        Then we see that $\Delta'' \in S$. If $\widetilde{G}\uind{i-1} \prec_{F\uind{n},\ldots, F\uind{i}} F\uind{i-1}$, then we replace $\Delta''$ with
        $$ \Delta''' = (F\uind{n}\supset \ldots \supset F\uind{i+1} \supset F\uind{i} \supset \widetilde{G}\uind{i-1} \supset \ldots \supset F\uind{i_{l}} \supset \ldots \supset F\uind{0}).$$
        Then we see that $\Delta''' \in S$ and this violates the minimality of $\Delta'$. On the contrary, if $\widetilde{G}\uind{i-1} \succeq_{F\uind{n},\ldots, F\uind{i}} F\uind{i-1}$, this means that
        $$ F\uind{i-1} \subset F\uind{i} \cap \bigcup_{H\uind{i} \prec_{F\uind{n},\ldots, F\uind{i+1}} F\uind{i}}H\uind{i}$$
        since $\widetilde{G}\uind{i-1} = F\uind{i} \cap \widetilde{G}\uind{i}$ is contained on the right hand side. However, this implies that $F\lsta\uind{i} \prec_{F\uind{n},\ldots, F\uind{i+1}} F\uind{i}$, which is a contradiction since $i_{1},\ldots, i_{l}$ are exactly those satisfying $F\uind{j}\lsta \prec_{F\uind{n},\ldots, F\uind{j+1}} F\uind{j}$ and $i$ is not one of them.

        The second case is when $i < i_{l}$. This implies that there exists $\widetilde{G}\uind{i} \prec_{F\uind{n},\ldots, F\uind{i+1}} F\uind{i}$ such that $\widetilde{G}\uind{i} \cap F\uind{i}$ is a facet of $F\uind{i}$. Then we follow the same lines as above and get a contradiction on the minimality of $i$.
    \end{proof}

    \subsection{The Decomposition Theorem} \label{section-Decompositiontheorem}
    In this subsection, we study the decomposition theorem for a proper birational toric morphism. We note that the singular cohomology of the fibers determine {\bf both} the intersection cohomology stalks and the coefficients in the decomposition theorem, and we describe how to compute both of these explicitly (see Remark \ref{rema:explicit-computation-of-H} below). We follow along the lines of \cite{CFS-Effectivedecompositiontheorem}, which discusses the decomposition theorem for Schubert varieties. However, the same argument works in the setting of toric varieties.
    
    Let $X$ be the affine toric variety of dimension $n$ associated to a full-dimensional rational polyhedral cone $\sigma$ and let $\Sigma_{X}$ be the associated fan. Consider a toric variety $Y$ obtained by a subdivision $\Sigma_{Y}$ of $\Sigma_{X}$, and let $\pi \colon Y \to X$ be the corresponding toric morphism. The decomposition theorem in this setting (see \cite[Theorem D]{dCMM-toricmaps}) tells us that
    \begin{equation}\label{equation:dec-thm-BBD}
        \bfR \pi\lsta \IC_{Y} \simeq \bigoplus_{\tau \in \Sigma_{X}} \bigoplus_{j\in \ZZ} \IC_{S_{\tau}}^{\oplus s_{\tau, j}}[-j].
    \end{equation}
    This is a more precise version of the decomposition theorem in \cite{BBD-Faisceauxpervers}, which is achieved by exploiting the action of the torus.

    For $\mu, \tau \in \Sigma_{X}$, we define
    \begin{align*}
        F_{\tau}(q) & = \sum_{j} h^{j}(\pi^{-1}(x_{\tau})) q^{j} \\
        H_{\mu, \tau}(q) & = \sum_{j} h^{j}(\IC_{S_{\mu}})_{x_{\tau}} q^{j} \\
        D_{\tau}(q) & = \sum_{j} s_{\tau, j}q^{j}.
    \end{align*}
    where $x_\tau \in O_\tau$ and
    $$h^{j}(\pi^{-1}(x_{\tau})) := \dim_{\QQ} H^{j}(\pi^{-1} (x_{\tau}), \QQ), \qquad \text{and} \qquad  h^{j}(\IC_{S_{\mu}})_{x_{\tau}} := \dim_{\QQ} \cH^{j}(\IC_{S_{\mu}})_{x_{\tau}}.$$
    Observe that $F_\tau$ is independent of the choice of $x_\tau$ by Proposition \ref{prop-fibertoric}, and $H_{\mu,\tau}$ is independent of the choice of $x_{\tau}$ by Remark \ref{rema:explicit-computation-of-H} below.
    
    We fix the notation for $F, H$, and $D$ from now on. If there is an ambiguity regarding the ambient toric variety, we will sometimes denote $H$ by $H^{X}$. We also point out that $H_{\mu, \tau}$ is nonzero only for $\mu\subset \tau$. If we additionally assume that $Y$ is a simplicial toric variety, we have $\IC_{Y} = \QQ_{Y}[n]$. By taking the stalk of both sides of Equation \ref{equation:dec-thm-BBD} at a point $x_{\tau} \in O_{\tau}$, we get
    \begin{equation} \label{equation-stalkdecomposition}
         F_{\tau}(q) q^{-n} = \sum_{\mu\subset \tau} H_{\mu, \tau}(q) \cdot D_{\mu}(q).
    \end{equation}
    Now we list some basic properties of these polynomials.
    \begin{enumerate}
        \item $D_{\tau}(q) = D_{\tau}(q^{-1})$ by Poincar\'e duality \cite{dCM-Decomposition}*{\S 1.6. (10)}.
        \item $q^{n- d_{\tau}}H_{\mu, \tau}(q)$ is strictly supported in negative degrees if $\mu \subsetneq \tau$ by \cite{dCM-Decomposition}*{\S 2.1.(12)}.
        \item $H_{\tau, \tau}(q) = q^{d_{\tau} - n}$.
        \item $D_{0}(q) = q^{0}$.
    \end{enumerate}

    For convenience, we put
    $$ \widetilde{F}_{\tau}(q) = q^{-d_{\tau}} F_{\tau}(q) \quad \text{and} \quad \widetilde{H}_{\mu, \tau}(q) = q^{n-d_{\tau}} H_{\mu, \tau}(q).$$

    \begin{lemm} \label{lemma-ICstalksforToricSubvar}
        For every pair of faces $\mu,\tau$ of $\sigma$ with $\mu \subset \tau$, we have the equality
        $$ \widetilde{H}_{\mu, \tau}^{X}(q) = \widetilde{H}_{0, \bar{\tau}}^{S_{\mu}}(q) = \widetilde{H}_{0, \bar{\tau}}^{V_{\mu,\tau}}(q) $$
        where $\bar{\tau}$ in the second term is considered as a face $\bar{\tau} \subset N_{\mathbb{R}}/ \langle\mu\rangle$ of the cone $\bar{\sigma} \subset N_{\mathbb{R}}/ \langle\mu\rangle$ and the $\bar{\tau}$ in the third term is considered as the cone $\bar{\tau} \subset \langle\tau\rangle / \langle\mu\rangle$.
    \end{lemm}
    \begin{proof}
        The first equality is straightforward. For the second equality, we use the description of $U_{\bar{\tau}} \subset S_{\mu}$. Note that $U_{\bar{\tau}} \simeq V_{\mu, \tau} \times O_{\tau}$. Since $\dim O_{\tau} = n - d_{\tau}$, we get
        $$ h^{j}(\IC_{S_{\mu}})_{x_{\tau}} = h^{j+\dim O_{\tau}} (\IC_{V_{\mu, \tau}})_{x_{\bar{\tau}}} $$
        due to the shift of cohomological degrees in the intersection complex. This proves the second equality.
    \end{proof}

    \begin{rema}\label{rema:explicit-computation-of-H}
        We describe how to explicitly compute $\widetilde{H}_{\mu, \tau}$ and $D_{\tau}$ in terms of the combinatorial data of the map $\pi$. Proposition \ref{prop-fibercohomology} tells us how to compute $F_{\tau}$ explicitly in terms of the combinatorial data of $\pi$. By induction on the dimension of $\tau$, we can assume that we have computed $D_{\mu}$ for all $\mu \subsetneq \tau$, since $D_{0}(q) = q^{0}$. Moreover, we can also assume that we have computed $\widetilde{H}_{\mu, \tau}$ for $0 \subsetneq \mu \subseteq \tau$ using Lemma \ref{lemma-ICstalksforToricSubvar} and by induction on the dimension of $\tau$. We proceed to compute $\widetilde{H}_{0, \tau}$ and $D_{\tau}$. By Equation (\ref{equation-stalkdecomposition}) and using $D_{0}(q) = q^{0}$ and $H_{\tau,\tau}(q) = q^{d_{\tau} - n}$, we get
    $$ F_{\tau}(q) q^{-n} = H_{0, \tau}(q) +  q^{d_{\tau} - n} D_{\tau}(q) + \sum_{0 \subsetneq \mu \subsetneq \tau} H_{\mu, \tau}(q) D_{\mu}(q). $$
    After multiplying by $q^{n-d_{\tau}}$, we get
    \begin{equation} \label{equation:stalk-decomposition-tilde}
        \widetilde{F}_{\tau}(q) - \sum_{0 \subsetneq \mu \subsetneq \tau} \widetilde{H}_{\mu, \tau}(q) D_{\mu}(q) = \widetilde{H}_{0, \tau}(q) + D_{\tau}(q).
    \end{equation}
    Note that $D_{\tau}(q)$ and $\widetilde{H}_{0, \tau}(q)$ can be completely determined provided that the left-hand side of Equation (\ref{equation:stalk-decomposition-tilde}) is known since $\widetilde{H}_{0, \tau}(q)$ is supported in strictly negative degrees and $D_{\tau}(q) = D_{\tau}(q^{-1})$. Hence, we can compute $D_{\tau}(q)$ and $\widetilde{H}_{0, \tau}(q)$ inductively. Additionally, we observe from the computation above that $\widetilde{H}_{\mu, \tau}$ does not depend on the choice of $x_{\tau}$.

    \noindent In the appendix, we use this idea to explicitly compute $\widetilde{H}_{\mu,\tau}$ in dimensions $\leq 4$.
    \end{rema}

    We finally upgrade the decomposition in (\ref{equation:dec-thm-BBD}) at the level of Hodge modules. Given a Hodge module $\cM$ of weight $w$ with filtration $F_\bullet$ and underlying perverse sheaf $K$, the \textit{Tate twist} of $M$ by an integer $k$ is a Hodge module, denoted $\cM(k)$, of weight $w-2k$ which has the same $\cD$-module structure as $\cM$ but with the filtration $F_{\bullet - k}$ and the perverse sheaf $K \otimes_{\QQ} (2\pi i)^k \QQ$. 

    \begin{prop}\label{prop:saito-dec-thm-toric}
        Let $X$ be the affine toric variety of dimension $n$ associated to a full-dimensional rational polyhedral cone $\sigma \subset N_{\RR}$ and let $\Sigma_{X}$ be the associated fan. Consider a toric variety $Y$ obtained by a \textbf{simplicial} subdivision $\Sigma_{Y}$ of $\Sigma_{X}$, and let $\pi \colon Y \to X$ be the corresponding toric morphism. Then, Saito's Decomposition theorem takes the form
        $$ \pi\lsta \IC_{Y}^{H} \simeq \pi\lsta \QQ_{Y}^{H}[n] \simeq \bigoplus_{j} \bigoplus_{\mu \subset \sigma} \left( \IC_{S_{\mu}}^{H} ( - \frac{d_{\mu} + j}{2}) \right)^{\oplus s_{\mu, j} } [-j].$$
    \end{prop}

    \begin{proof}        
        Saito's Decomposition theorem tells us that
        \begin{equation}\label{equation:saito's-Dec-Thm-preliminary}
            \pi\lsta \QQ^H_{Y}[n] \simeq \bigoplus_{j} \bigoplus_{\mu \subset \sigma} M_{\mu, j} [-j]
        \end{equation}
        where $M_{\mu,j}$ is a Hodge module on $X$ supported on $S_\mu$, of weight $n+j$, with its underlying perverse sheaf equal to $\IC_{S_\mu}^{\oplus s_{\mu,j}}$. We need to show that $M_{\mu,j} \simeq \left( \IC_{S_{\mu}}^{H} ( - \frac{d_{\mu} + j}{2}) \right)^{\oplus s_{\mu, j} }$.

        \textbf{Step 1.} We first do the case $\mu = \sigma$. In this case, we denote $S_\sigma$ by $x_\sigma$ since $S_\sigma$ is just the torus fixed point of $X$. Observe that $M_{\sigma,j}$ is a Hodge structure with the underlying vector space given by $\QQ_{x_\sigma}^{\oplus s_{\sigma,j}}$. We need to show that $M_{\sigma,j}$ is of Hodge--Tate type, i.e. it is isomorphic to the Hodge structure $\left( \QQ_{x_\sigma}^{H}(-\frac{n+j}{2})\right)^{\oplus s_{\sigma, j}}$. We have the diagram
        \[
        \begin{tikzcd}
            E \ar[r, "i_E"] \ar[d, swap, "\pi_E"] & Y \ar[d, "\pi"]\\
            x_{\sigma} \ar[r, "i_{x_\sigma}"] & X
        \end{tikzcd}
        \]
        where $E$ is the reduced inverse image of $x_\sigma$. Consider the pullback of both sides of (\ref{equation:saito's-Dec-Thm-preliminary}) to $x_\sigma$
        \begin{equation}\label{equation:saito's-Dec-Thm-preliminary-restriction}
            i^*_{x_\sigma}\pi\lsta \QQ^H_{Y}[n] \simeq i_{x_\sigma}^* \left( \bigoplus_{j} \bigoplus_{\mu \subset \sigma} M_{\mu, j} [-j] \right).
        \end{equation}
        By proper base change for Hodge modules (see \cite{saito1990mixedHodgemodules}*{4.4.3}), the left hand side becomes
        \[ i_{x_\sigma}^*(\pi\lsta \QQ_{Y}^{H}[n]) \simeq (\pi_E)\lsta (i_E^*\QQ^H_Y[n]) \simeq (\pi_E)\lsta \QQ^H_E[n]  \]
        whose cohomologies are the singular cohomology groups $H^p(E,\QQ)$, while the right hand side has $M_{\sigma,j}[-j]$ as one of the direct summands since $i_{x_\sigma}^*M_{\sigma,j} = M_{\sigma,j}$. Now, the singular cohomology groups $H^p(E,\QQ)$ are pure of weight $p$, and of Hodge--Tate type (see Proposition \ref{prop-fibercohomology}). Taking the $j$-th cohomology on both sides of (\ref{equation:saito's-Dec-Thm-preliminary-restriction}), we see that $M_{\sigma,j}$ is a direct summand of $H^{n+j}(E,\QQ)$, and so we conclude that $M_{\sigma,j}$ is pure of weight $n+j$, and of Hodge--Tate type, i.e. $M_{\sigma,j} \simeq \left( \QQ_{x_\sigma}^{H}(-\frac{n+j}{2})\right)^{\oplus s_{\sigma, j}}$ as required.

        \textbf{Step 2.} Now for the general case, consider $M_{\mu_0,j}$ for a face ${\mu_0} \subset \sigma$. Denote the affine open subset of $X$ corresponding to ${\mu_0}$ by $U_{\mu_0}$. Observe that $U_{\mu_0} \cap S_{\mu_0} = O_{\mu_0}$, the orbit in $X$ corresponding to ${\mu_0}$. Thus, it suffices to show that 
        $$M_{{\mu_0},j}|_{U_{\mu_0}} \simeq \left( \QQ_{O_{{\mu_0}}}^{H} ( - \frac{j}{2}) \right)^{\oplus s_{{\mu_0}, j} }$$
        since we know that $M_{{\mu_0},j}$ is the intermediate extension of its restriction to $U_{\mu_0}$. Just as in Remark \ref{rema:Basics-Toric-varieties}(2), we fix a non-canonical splitting $N = N_{\mu_0} \oplus (N \cap \langle {\mu_0} \rangle)$ to get $U_{\mu_0} \simeq V_{\mu_0} \times O_{\mu_0}$ where $V_{\mu_0}$ is the toric variety corresponding to the cone ${\mu_0} \subset (N \cap \langle {\mu_0} \rangle)_{\RR}$. Denote the inverse image of $U_{\mu_0}$ under $\pi$ as $Y_{\mu_0} = \pi^{-1}(U_{\mu_0})$. By using the same non-canonical splitting, we get that $Y_{\mu_0} \simeq W_{\mu_0} \times O_{\mu_0}$, where $p:W_{\mu_0} \to V_{\mu_0}$ is the induced map and the following diagram commutes
        \[
        \begin{tikzcd}
            Y_{\mu_0} \ar[r, "\pr"] \ar[d, "\pi"] & W_{\mu_0} \ar[d, "p"]\\
            U_{\mu_0} \ar[r, "\pr"] & V_{\mu_0}
        \end{tikzcd}
        \]
        where $\pr$ denotes the projection onto the first factor. Observe now that $\pr^*\QQ^H_{W_{\mu_0}} = \QQ^H_{Y_{\mu_0}}$. Thus we get
        \[ \pi_*\QQ^H_{Y_{\mu_0}}[n] \simeq \pi_*\left( \pr^*\QQ^H_{W_{\mu_0}}[n] \right) \simeq \pr^*\left(p_*\QQ^H_{W_{\mu_0}}[n] \right). \]
        We apply Saito's decomposition theorem to the map $p$ to get
        \[ p\lsta \QQ^H_{W_{\mu_0}}[d_{\mu_0}] \simeq \bigoplus_{j} \bigoplus_{\tau \subset {\mu_0}} M^{V_{\mu_0}}_{\tau, j} [-j]. \]
        Combining the above two equations, we get
        \begin{align*}
            \pi_*\QQ^H_{Y_{\mu_0}}[n] &\simeq \pr^* \left( \bigoplus_{j} \bigoplus_{\tau \subset {\mu_0}} M^{V_{\mu_0}}_{\tau, j} [-j] \right)[n-d_{\mu_0}]\\
            &\simeq \bigoplus_{j} \bigoplus_{\tau \subset {\mu_0}} \pr^*(M^{V_{\mu_0}}_{\tau, j}) [-j+n-d_{\mu_0}].
        \end{align*}
        On the other hand, restricting both sides of (\ref{equation:saito's-Dec-Thm-preliminary}) to $U_{\mu_0}$ gives
        \[ \pi_*\QQ^H_{Y_{\mu_0}}[n] \simeq (\pi\lsta \QQ^H_{Y}[n])|_{U_{{\mu_0}}} \simeq \bigoplus_{j} \bigoplus_{\tau \subset \sigma} (M_{\tau, j})|_{U_{\mu_0}} [-j] \simeq \bigoplus_{j} \bigoplus_{{\tau} \subset \mu_0} (M_{{\tau}, j})|_{U_{\mu_0}} [-j] \]
        where the last isomorphism comes from the fact that $U_{\mu_0} \cap S_\tau = \emptyset$ for $\tau \not\subset \mu_0$. By comparing the above two decompositions of $\pi_*\QQ^H_{Y_{\mu_0}}[n]$, we get that
        \[ (M_{{\tau}, j})|_{U_{\mu_0}} \simeq \pr^*(M^{V_{\mu_0}}_{\tau, j-n+d_{\mu_0}}) \text{ for all } \tau \subset \mu_0 . \]
        In particular, for $\tau = \mu_0$, we get that
        \[ (M_{{\mu_0}, j})|_{U_{\mu_0}} \simeq \pr^*(M^{V_{\mu_0}}_{\mu_0, j-n+d_{\mu_0}}). \]
        Step 1 applied to the map $p:W_{\mu_0} \to V_{\mu_0}$ implies that the factor $M^{V_{\mu_0}}_{{\mu_0},j-n+d_{\mu_0}}$ which lives over the torus fixed point $x_{\mu_0}$ of $V_{\mu_0}$, is isomorphic to the Hodge structure $\left( \QQ_{x_{\mu_0}}^{H}(-\frac{j}{2}) \right)^{\oplus s^{V_{\mu_0}}_{\tau,j-n+d_{\mu_0}}}$. Therefore, we get $s^{V_{\mu_0}}_{\tau,j-n+d_{\mu_0}} = s_{\mu_0,j}$ and
        \[ (M_{{\mu_0}, j})|_{U_{\mu_0}} \simeq \pr^*\left( \QQ_{x_{\mu_0}}^{H}(-\frac{j}{2}) \right)^{\oplus s^{V_{\mu_0}}_{\tau,j-n+d_{\mu_0}}} \simeq \left( \QQ_{O_{{\mu_0}}}^{H} ( - \frac{j}{2}) \right)^{\oplus s_{{\mu_0}, j} }, \]
        where the last isomorphism follows from the fact that $\pr^{-1}(x_{\mu_0}) = O_{\mu_0} \subset U_{\mu_0}$.
    \end{proof}

    \section{Higher direct images of Kähler differentials for toric morphisms} \label{section:Kahler-diff}

    Let $\pi \colon Y\to X$ be a barycentric resolution of an affine toric variety $X$ as in Remark \ref{rema:barycentric-resolution}. The main goal of this section is to describe a systematic way to compute the higher direct images of the sheaves of reflexive differentials.

    The Ishida complex (\ref{equation:Ishida-complex}) provides a resolution of $\Omega\ubind{p}_Y$.
    The terms appearing in the Ishida complex are structure sheaves of various torus invariant subvarieties. Combining this fact with Remark \ref{rema:Basics-Toric-varieties} (1), we see that the higher direct images of $\Omega_{Y}\ubind{p}$ can be computed by calculating the cohomology of the pushforward of the Ishida complex along $\pi$. The pushforward of the Ishida complex is given by:
    \begin{equation}\label{equation:pushforward-of-Ishida-complex}
            0 \to \Omega_{X}\ubind{p} \to V_{0}^{p} \otimes \cO_{X} \to \bigoplus_{\substack{\nu \in \Sigma_{Y} \\ \dim \nu = 1}} V_{\nu}^{p} \otimes \cO_{S_{\pi_*(\nu)}} \to \bigoplus_{\substack{\nu \in \Sigma_{Y} \\ \dim \nu = 2}} V_{\nu}^{p} \otimes \cO_{S_{\pi_*(\nu)}} \to \ldots \to \bigoplus_{\substack{\nu \in \Sigma_{Y} \\ \dim \nu = p}} V_{\nu}^{p} \otimes \cO_{S_{\pi_*(\nu)}}\to 0.
    \end{equation}

    We first describe the morphisms in the complex. Let $\mu , \nu \in \Sigma_{Y}$ such that $\dim \mu = l$ and $\dim \nu = l+1$. Then we see that $S_{\pi\lsta \nu} \subset S_{\pi\lsta \mu}$ if and only if $\mu \subset \nu$. Let $\rho_{1},\ldots, \rho_{l}$ be the rays of $\mu$. If $\mu \subset \nu$, there is a unique ray $\rho_{l+1}$ such that $\nu$ is the span of $\rho_{1},\ldots, \rho_{l},\rho_{l+1}$. Recall from Remark \ref{rema:Basics-Toric-varieties} (4) that
    $$ V_{\mu}^{p} = \left( \bigwedge^{p-l} \mu^{\perp} \right) \otimes \frac{M_{\CC}}{\rho_{1}^{\perp}} \otimes \ldots \otimes \frac{M_{\CC}}{\rho_{l}^{\perp}}.$$
    Note that $\mu^{\perp} / \nu^{\perp}$ can be naturally identified with $M_{\CC} / \rho_{l+1}^{\perp}$. From the short exact sequence $0 \to \nu^{\perp} \to \mu^{\perp} \to \mu^{\perp} / \nu^{\perp} \to 0$, we get a surjection
    $$ \bigwedge^{p-l} \mu^{\perp} \to \bigwedge^{p-l-1} \nu^{\perp} \otimes \frac{\mu^{\perp}}{\nu^{\perp}} \simeq \bigwedge^{p-l-1} \nu^{\perp} \otimes \frac{M_{\CC}}{\rho_{l+1}^{\perp}}. $$
    The map between $V_{\mu}^{p} \otimes \cO_{S_{\pi\lsta \mu}} \to V_{\nu}^{p} \otimes \cO_{S_{\pi\lsta \nu}}$ in (\ref{equation:pushforward-of-Ishida-complex}) is given by the map above (up to a well-defined sign after choosing the order of the rays) tensored with $\frac{M_{\CC}}{\rho_{1}^{\perp}} \otimes \ldots \otimes \frac{M_{\CC}}{\rho_{l}^{\perp}}$ and the restriction map $\cO_{S_{\pi\lsta \mu}} \to \cO_{S_{\pi\lsta \nu}}$ if $\mu \subset \nu$. Otherwise, the map is zero.
    
    Using the action of the torus, the complex in (\ref{equation:pushforward-of-Ishida-complex}) decomposes into various eigenspaces, hence it carries a natural $M$-grading. For example, for $u \in M$ viewed as a character of the torus, the degree $u$-part of $\Omega_{X}\ubind{p}$ can be described as
    $$ \left(\Omega_{X}\ubind{p}\right)_{u} = \{ \alpha \in \Omega_{X}\ubind{p} : g\sta \alpha = u(g)\cdot \alpha \text{ for all } g \in \Spec \CC[M] \}\\ $$
    after identifying $\Omega_{X}\ubind{p}$ with its space of global sections. First, let us consider this complex in degree $0 \in M$. Since every $\cO_{S_{\pi_*(\nu)}}$ contains $\chi^0$, the pushforward of the Ishida complex in degree 0 is
    \begin{equation}\label{equation:pushforward-of-Ishida-complex-grading-0}
            0 \to V_{0}^{p} \to \bigoplus_{\substack{\nu \in \Sigma_{Y} \\ \dim \nu = 1}} V_{\nu}^{p} \to \bigoplus_{\substack{\nu \in \Sigma_{Y} \\ \dim \nu = 2}} V_{\nu}^{p} \to \ldots \to \bigoplus_{\substack{\nu \in \Sigma_{Y} \\ \dim \nu = p}} V_{\nu}^{p} \to 0,
        \end{equation}
    with $V^p_0$ in cohomological degree $0$.

    Let us write this in terms of the polytopal complex $\cC_Y$ instead of $\Sigma_Y$. For $\gamma \in \cC_{Y}$, we denote by $V_{\gamma}^{p}$ the vector space $V_{\mathrm{span}_{\RR_{\geq 0}}\gamma}^{p}$ where $\mathrm{span}_{\RR_{\geq 0}} \gamma$ is the cone in $\Sigma_{Y}$ corresponding to $\gamma$.  We point out that there is a difference of dimension by 1 between the dimension as a cone and the dimension as a simplex, hence (\ref{equation:pushforward-of-Ishida-complex-grading-0}) becomes
    \begin{equation}\label{equation:pushforward-of-Ishida-complex-grading-0-polytope}
             0 \to V_{0}^{p} \to \bigoplus_{\substack{\gamma \in \cC_{Y} \\ \dim \gamma = 0}} V_{\gamma}^{p} \to \bigoplus_{\substack{\gamma \in \cC_{Y} \\ \dim \gamma = 1}} V_{\gamma}^{p} \to \ldots \to \bigoplus_{\substack{\gamma \in \cC_{Y} \\ \dim \gamma = p-1}} V_{\gamma}^{p} \to 0.
        \end{equation}

    We choose the lexicographic shelling of $\cC_Y$ as described in Proposition \ref{prop:barycentric-subdivision-shellable}. Let $\Delta_1,\dots,\Delta_m$ denote the shelling. 

    \begin{rema}\label{remark:exact-complex-of-simplex}
        Let us describe the complex associated to a simplex of type $l$ for $l \leq p$. Let $\Delta_k$ be a simplex of type $l$ with $l\leq p$. Let $\alpha = [v\lind{i_{1}},\ldots, v\lind{i_{l}}]$ denote the face of $\Delta_{k}$ that is not a face of one of $\Delta_{1},\ldots, \Delta_{k-1}$. We associate a Koszul-type complex to $\alpha$:
        \begin{equation}\label{equation:Exact-diagram-quotient-off}
           0 \to V_{\alpha}^{p} \to \bigoplus_{\substack{\gamma \in \Delta_k \\ \dim \gamma = l \\ \alpha \subset \gamma }} V_{\gamma}^{p} \to \bigoplus_{\substack{\gamma \in \Delta_k \\ \dim \gamma = l+1 \\ \alpha \subset \gamma }} V_{\gamma}^{p} \to \ldots \to \bigoplus_{\substack{\gamma \in \Delta_k \\ \dim \gamma = p-1 \\ \alpha \subset \gamma }} V_{\gamma}^{p} \to 0  
        \end{equation}
        with $V^p_\alpha$ in cohomological degree $l$. We note that the morphisms in the complex come from the push-forward of the Ishida complex in Equation \ref{equation:pushforward-of-Ishida-complex}. When $l=p$, this is a complex concentrated in a single degree, namely $V^p_\alpha[-p]$. For $l< p$, this is exactly \cite{toric-SVV}*{Equation (3)} tensored by $\frac{V}{v\lind{i_{1}}^{\perp}} \otimes \ldots \otimes \frac{V}{v\lind{i_{l}}^{\perp}}$, and is hence exact. In either case, observe that the complex is exact in all cohomological degrees other than $p$.
    \end{rema}

    \begin{prop}\label{prop:pushforward-of-Ishida-complex-grading-0}
        The complex (\ref{equation:pushforward-of-Ishida-complex-grading-0-polytope}) is exact in all cohomological degrees other than $p$.
    \end{prop}

    \begin{proof}
        Let $\Delta_{1},\ldots, \Delta_{m}$ be a shelling of $\cC_{Y}$. We compute the cohomology of the complex (\ref{equation:pushforward-of-Ishida-complex-grading-0-polytope})
        $$ A^{\bullet} := 0 \to V_{0}^{p} \to \bigoplus_{\substack{\gamma \in \cC_{Y} \\ \dim \gamma = 0}} V_{\gamma}^{p} \to \bigoplus_{\substack{\gamma \in \cC_{Y} \\ \dim \gamma = 1}} V_{\gamma}^{p} \to \ldots \to \bigoplus_{\substack{\gamma \in \cC_{Y} \\ \dim \gamma = p-1}} V_{\gamma}^{p} \to 0,$$
        with $V_{0}^{p}$ in cohomological degree $0$, by putting a filtration on $A^{\bullet}$ from the shelling and computing the associated spectral sequence. We define $F^{i}A^{\bullet}$ by
        $$ F^{i}A^{\bullet}:= 0 \to 0 \to \bigoplus_{\substack{\gamma \in \cC_{Y}\setminus \bigcup_{k=1}^{i}\Delta_k \\ \dim \gamma = 0}} V_{\gamma}^{p} \to \bigoplus_{\substack{\gamma \in \cC_{Y}\setminus \bigcup_{k=1}^{i}\Delta_k \\ \dim \gamma = 1}} V_{\gamma}^{p} \to \ldots \to \bigoplus_{\substack{\gamma \in \cC_{Y}\setminus \bigcup_{k=1}^{i}\Delta_k \\ \dim \gamma = p-1}} V_{\gamma}^{p} \to 0$$
        for $i > 0$, and $F^{0}A^{\bullet} = A^{\bullet}$. First we observe that $F^{i}A^{\bullet}$ are indeed subcomplexes of $A^{\bullet}$. For this, it is enough to show that if $\gamma \in \cC_{Y} \setminus \bigcup_{k=1}^{i} \Delta_{k}$ with $\dim \gamma = l$, and if $\gamma' \in \bigcup_{k=1}^{i} \Delta_{k}$ with $\dim \gamma' = l+1$, then the morphism $V_{\gamma}^{p} \to V_{\gamma'}^{p}$ in $A^{\bullet}$ is zero. But this is true because $\gamma$ cannot be a face of $\gamma'$ in this situation, and hence $F^{i}A^{\bullet}$ is indeed a subcomplex of $A^{\bullet}$. We describe the graded pieces of this filtration. First, notice that
        \begin{align*}
           \gr_{F}^{0}A^{\bullet} & = 0 \to V_{0}^{p} \to \bigoplus_{\substack{\gamma \in \Delta_1 \\ \dim \gamma = 0}} V_{\gamma}^{p} \to \bigoplus_{\substack{\gamma \in \Delta_1 \\ \dim \gamma = 1}} V_{\gamma}^{p} \to \ldots \to \bigoplus_{\substack{\gamma \in \Delta_1 \\ \dim \gamma = p-1}} V_{\gamma}^{p} \to 0.
        \end{align*}
        For $i> 0$, let $\alpha = [v\lind{i_{1}}, \ldots, v\lind{i_{l}}]$ be the face of $\Delta_{i}$ that is not a face of any of $\Delta_{1},\ldots, \Delta_{i-1}$. Notice that
        $$ \gr_{F}^{i} A^{\bullet} = 0 \to V_{\alpha}^{p} \to \bigoplus_{\substack{\gamma \in \Delta_k \setminus \bigcup_{i=1}^{k-1} \Delta_i \\ \dim \gamma = l}} V_{\gamma}^{p} \to \bigoplus_{\substack{\gamma \in \Delta_k \setminus \bigcup_{i=1}^{k-1} \Delta_i \\ \dim \gamma = l+1}} V_{\gamma}^{p} \to \ldots \to \bigoplus_{\substack{\gamma \in \Delta_k \setminus \bigcup_{i=1}^{k-1} \Delta_i \\ \dim \gamma = p-1}} V_{\gamma}^{p} \to 0,$$
        where $V_{\alpha}^{p}$ sits in cohomological degree $l$. Observe that this is exactly the complex described in Remark \ref{remark:exact-complex-of-simplex}. Consider the spectral sequence
        $$ E_{1}^{i,j} = H^{i+j} \gr_{F}^{i} A^{\bullet} \implies H^{i+j} A^{\bullet}.$$
        We point out that $H^{i+j} \gr_{F}^{i} A^{\bullet} = 0$ unless $i + j = p$ by Remark \ref{remark:exact-complex-of-simplex}. Therefore, the spectral sequence degenerates at $E_{1}$, and the only non-trivial cohomology of $A^{\bullet}$ is in cohomological degree $p$.
    \end{proof}

    \begin{rema}\label{rema:vector-spaces-barycentric-subdivision}
        Observe that Proposition \ref{prop:pushforward-of-Ishida-complex-grading-0} is actually a statement about a complex of vector spaces associated to the barycentric subdivision of the cone $\sigma$. In particular, this statement also holds for any face $\tau$ of $\sigma$ and its barycentric subdivision as well.
    \end{rema}

    Finally, let us see what happens in a degree $u \in \tau_{\circ}^*$. Observe that the terms $\cO_{S_{\pi_*(\nu)}}$ that are non-zero in degree $u$ are precisely those for which $\nu \subset \tau$. Therefore the pushforward of the Ishida complex in degree $u$ is
    \begin{equation}\label{equation:pushforward-of-Ishida-complex-grading-u}
            0 \to V_{0}^{p} \to \bigoplus_{\substack{\nu \in \Sigma_{Y} \\ \nu \subset \tau,\dim \nu = 1}} V_{\nu}^{p} \to \bigoplus_{\substack{\nu \in \Sigma_{Y} \\ \nu \subset \tau,\dim \nu = 2}} V_{\nu}^{p} \to \ldots \to \bigoplus_{\substack{\nu \in \Sigma_{Y} \\ \nu \subset \tau,\dim \nu = p}} V_{\nu}^{p} \to 0.
        \end{equation} 
    Let $W=\tau^\perp$ and let $\overline{V} = V/\tau^\perp$ (note that $\dim W = n-d_\tau$ and $\dim \overline{V} = d_{\tau}$). Fix a non-canonical splitting $V = W \oplus \overline{V}$. Observe that we have a non-canonical isomorphism
    \[ \bigwedge^l V \simeq \bigoplus_{i=0}^{d_\tau} \left(\bigwedge^{l-i}W \otimes \bigwedge^i \overline{V} \right), \]
    with the convention that $\displaystyle \bigwedge^j W =0$ for $j>n-d_\tau$.

    More generally, given $\nu \subset \tau$, observe that $W = \tau^\perp \subset \nu^\perp$. Define $\overline{V}_\nu := \nu^\perp/\tau^\perp$. The non-canonical splitting of $V$ we fixed earlier induces a splitting $\nu^\perp = W \oplus \overline{V}_\nu$. We can now write the pushforward of the Ishida complex in degree $u$ as a direct sum
    \begin{align*}
            0 \to &V_{0}^{p} \to \bigoplus_{\substack{\nu \in \Sigma_{Y} \\ \nu \subset \tau,\dim \nu = 1}} V_{\nu}^{p} \to \bigoplus_{\substack{\nu \in \Sigma_{Y} \\ \nu \subset \tau,\dim \nu = 2}} V_{\nu}^{p} \to \ldots \to \bigoplus_{\substack{\nu \in \Sigma_{Y} \\ \nu \subset \tau,\dim \nu = p}} V_{\nu}^{p} \to 0\\
            & \simeq \bigoplus_{i=0}^{d_\tau} \left(\bigwedge^{p-i}W \otimes \left( \bigwedge^{i} \overline{V} \to \bigoplus_{\substack{\nu \in \Sigma_{Y} \\ \nu \subset \tau, \dim \nu = 1}} \overline{V}_{\nu}^{i} \to \bigoplus_{\substack{\nu \in \Sigma_{Y} \\ \nu \subset \tau, \dim \nu = 2}} \overline{V}_{\nu}^{i} \to \ldots \to \bigoplus_{\substack{\nu \in \Sigma_{Y} \\ \nu \subset \tau, \dim \nu = i}} \overline{V}_{\nu}^{i} \right) \right)
        \end{align*} 

    Now, observe that the complex 
    \[ \bigwedge^{i} \overline{V} \to \bigoplus_{\substack{\nu \in \Sigma_{Y} \\ \nu \subset \tau, \dim \nu = 1}} \overline{V}_{\nu}^{i} \to \bigoplus_{\substack{\nu \in \Sigma_{Y} \\ \nu \subset \tau, \dim \nu = 2}} \overline{V}_{\nu}^{i} \to \ldots \to \bigoplus_{\substack{\nu \in \Sigma_{Y} \\ \nu \subset \tau, \dim \nu = i}} \overline{V}_{\nu}^{i} \]
    is the complex of vector spaces associated to the barycentric subdivision of the cone $\tau$. Therefore by Remark \ref{rema:vector-spaces-barycentric-subdivision}, it follows that this complex is exact at all cohomological degrees other than $i$.
    
    Now, we define the generating function for the pushforward of the sheaves of reflexive differentials. For $\tau \in \Sigma_{X}$, we define
    $$ \Omega_{\tau}(K, L) = \sum_{k, l} \dim_{\CC} (R^{n+k+l}\pi\lsta \Omega_{Y}\ubind{-k})_{u} K^{k} L^{l},$$
    where $u \in \tau_{\circ}\sta$. This is independent of the choice of $u \in \tau_{\circ}\sta$, as shown by the following combinatorial formula for $\Omega_{\tau}(K, L)$.

    \begin{prop} \label{prop-formulaKahlerdiff}
        If $\pi \colon Y \to X$ is a barycentric resolution of $X$, then
        $$ \Omega_{\tau}(K, L) = L^{-n}(1 + K^{-1}L)^{n- d_{\tau}} \sum_{\substack{\mu \subset \tau \\ \mu \in \Sigma_{X}}} \left( \sum_{j = 0}^{d_{\tau}} d_{j}(\mu)  (1-K^{-1}L^{2})^{d_{\tau} - j} (K^{-1}L^{2})^{j} \right). $$
        Moreover, $\Omega_{\tau}$ is related to $F_{\tau}$ by the following formula:
        $$ \Omega_{\tau}(K, L) = L^{-n} (1 + K^{-1}L)^{n- d_{\tau}} F_{\tau}(LK^{- \frac{1}{2}}). $$
    \end{prop}

    \begin{proof}
        We compute the dimension of $(R^{i} \pi\lsta \Omega_{Y}\ubind{p})_{u}$. Let $W = \tau^{\perp}$ and $\overline{V} = V / \tau^{\perp}$. Note that by Proposition \ref{prop:pushforward-of-Ishida-complex-grading-0}, the cohomology of the complex
        $$ \bigwedge^{p-i}W \otimes \left( \bigwedge^{i} \overline{V} \to \bigoplus_{\substack{\nu \in \Sigma_{Y} \\ \nu \subset \tau, \dim \nu = 1}} \overline{V}_{\nu}^{i} \to \bigoplus_{\substack{\nu \in \Sigma_{Y} \\ \nu \subset \tau, \dim \nu = 2}} \overline{V}_{\nu}^{i} \to \ldots \to \bigoplus_{\substack{\nu \in \Sigma_{Y} \\ \nu \subset \tau, \dim \nu = i}} \overline{V}_{\nu}^{i} \right) $$
        is concentrated in degree $i$ and is equal to $(R^{i}\pi\lsta \Omega_{Y}\ubind{p})_{u}$. Furthermore, the number of summands of the $j$-th term is equal to $\sum_{\mu \subset \tau, \mu \in \Sigma_{X}} d_{j}(\mu)$ and the dimension of each summand of the $j$-th term is equal to ${d_{\tau} -j\choose i-j}$. Hence we have
        \begin{align*}
            & \dim_{\CC} (R^{i}\pi\lsta \Omega_{Y}\ubind{p})_{u}\\& = {n - d_{\tau} \choose p-i} \cdot \left( (-1)^{i} {d_{\tau} \choose i} \sum_{\substack{\mu \subset \tau \\ \mu \in \Sigma_{X}}} d_{0}(\mu)  +(-1)^{i-1} {d_{\tau} -1 \choose i-1} \sum_{\substack{\mu \subset \tau \\ \mu \in \Sigma_{X}}} d_{1} (\mu) +\cdots  + {d_{\tau} - i \choose 0} \sum_{\substack{\mu \subset \tau \\ \mu \in \Sigma_{X}}} d_{i}(\mu) \right) \\
            & = {n -d_{\tau} \choose p -i} \sum_{j = 0}^{i} \sum_{\substack{\mu \subset \tau \\ \mu \in \Sigma_{X}}} (-1)^{i-j} {d_{\tau} - j \choose i - j} d_{j}(\mu).
        \end{align*}

        Therefore, we have
        \begin{align*}
            \Omega_{\tau}(K, L) & = \sum_{k, l} \dim_{\CC} (R^{n+k+l} \pi\lsta \Omega_{Y}\ubind{-k})_{u} K^{k} L^{l} \\
            & = \sum_{k, l, j} \sum_{\mu \subset \tau} {n- d_{\tau}\choose -n-2k-l} (-1)^{n+k+l-j} 
            {d_{\tau} -j \choose n+k+l-j} d_{j}(\mu) K^{k} L^{l} \\
            & = \sum_{c, j, \mu \subset \tau} {n - d_{\tau} \choose -n-c} d_{j}(\mu) \sum_{k} (-1)^{n+c - j - k} {d_{\tau} - j \choose n+c-j-k} K^{k} L^{c- 2k}\\
            &\qquad (\text{where we set }c = 2k + l)\\
            & = \sum_{c, j, \mu \subset \tau} {n - d_{\tau} \choose -n-c} d_{j}(\mu) \left(\sum_{k'} (-1)^{k'} {d_{\tau} - j \choose k'} (K^{-1}L^{2})^{k'} \right) K^{n+c - j} L^{-2n-c + 2j} \\
            & \qquad (\text{where we set } k' = n+c -j-k) \\
            & = \sum_{c, j, \mu \subset \tau} {n - d_{\tau} \choose -n-c} d_{j}(\mu) (1 - K^{-1}L^{2})^{d_{\tau} - j} K^{n+c - j} L^{-2n-c + 2j}\\
            & =  L^{-n} \sum_{\mu \subset \tau} \sum_{j} d_{j}(\mu) (1 - K^{-1}L^{2})^{d_{\tau} - j} (K^{-1}L^{2})^{j}\sum_{c'} {n - d_{\tau} \choose c'} (K^{-1}L)^{c'}\\
            &\qquad (\text{where we set }c' = -n-c)\\
            & = L^{-n}(1 + K^{-1}L)^{n- d_{\tau}} \sum_{\substack{\mu \subset \tau \\ \mu \in \Sigma_{X}}}  \sum_{j}d_{j}(\mu) (1-K^{-1}L^{2})^{d_{\tau} - j} (K^{-1}L^{2})^{j} .
        \end{align*}
    
    This concludes the proof of the first part of the proposition. Now, we prove the second statement. By construction of the barycentric subdivision, we have
        $$ d_{j}(\tau) = \sum_{\mu \subsetneq \tau} d_{j-1}(\tau).$$
        If $t^{2} = K^{-1}L^{2}$, then we get
        \begin{align*}
            & \sum_{\mu \subset \tau}  \sum_{j} d_{j}(\mu) (1 - t^{2})^{d_{\tau} - j} (t^{2})^{j} \\
            & = \sum_{j} (d_{j}(\tau) + d_{j+1}(\tau)) (1-t^{2})^{d_{\tau} - j}(t^{2})^{j} \\
            & = \sum_{j} d_{j}(\tau) \left( (1-t^{2})^{d_{\tau}-j} (t^{2})^{j} + (1 - t^{2})^{d_{\tau} - j+1}(t^{2})^{j-1} \right) \\
            & = \sum_{j} d_{j}(\tau) (t^{2})^{j-1} (1 - t^{2})^{d_{\tau} - j}.
        \end{align*}
        Observe that the fiber $\pi^{-1}(x_{\tau})$ is an irreducible simplicial toric variety of dimension $d_{\tau} - 1$. First, we show that $\pi^{-1}(O_{\tau})$ is irreducible. Let $\rho_\tau$ denote the unique ray in $\Sigma_Y$ contained in the relative interior of $\tau$, then we have
        $$\pi^{-1}(O_{\tau}) = \displaystyle \bigsqcup_{\substack{\pi_*\mu = \tau \\ \mu \in \Sigma_{Y}}} O_\mu \subset \bigsqcup_{\substack{\rho_\tau \subset \mu \\ \mu \in \Sigma_{Y}}} O_\mu = \overline{O_{\rho_\tau}}.$$
        Since $\overline{O_{\rho_\tau}}$ is irreducible, and since $\pi^{-1}(O_{\tau}) \subset \overline{O_{\rho_\tau}}$ is an open subset (it is the inverse image of the open set $O_\tau$ under the natural map $\overline{O_{\rho_\tau}} \to \overline{O_\tau})$, we have that $\pi^{-1}(O_{\tau})$ is irreducible. By Proposition \ref{prop-fibertoric}, we have $\pi^{-1}(O_{\tau}) \simeq \pi^{-1}(x_{\tau}) \times O_{\tau}$. Hence $\pi^{-1}(x_{\tau})$ is irreducible and hence, irreducible simplicial by Proposition \ref{prop-fibertoric}.

        Now by Proposition \ref{prop-fibercohomology}, we have:
        $$ F_{\tau}(q) = \sum_{j} d_{j}(\tau) (q^{2}-1)^{d_{\tau} -j} .$$
        Since a proper simplicial toric variety satisfies Poincaré duality, we know that $F_{\tau}(q) = F_{\tau}(q^{-1}) q^{2 (d_{\tau} - 1)}$, and therefore
        $$ \sum_{j} d_{j}(\tau) (q^{-2} - 1)^{d_{\tau} - j} q^{2(d_{\tau} - 1)} = \sum_{j} d_{j}(\tau) (1 - q^{2})^{d_{\tau}-j}(q^{2})^{j-1}.$$
        Hence, we get
        $$ \sum_{\substack{\mu \subset \tau \\ \mu \in \Sigma_{X}}} \sum_{j} d_{j}(\mu) (1-K^{-1}L^{2})^{d_{\tau} - j} (K^{-1}L^{2})^{j} = \sum_{j} d_{j}(\tau) (K^{-1}L^{2} - 1)^{d_{\tau} - j}.$$
        Therefore,
        $$ \Omega_{\tau}(K,L) = L^{-n} (1+K^{-1}L)^{n- d_{\tau}} \sum_{j} d_{j}(\tau) (K^{-1}L^{2}-1)^{d_{\tau} - j}.$$
    \end{proof}

    \begin{rema}
        Even though Proposition \ref{prop-formulaKahlerdiff} is stated for barycentric resolutions, this method provides a general framework of computing the higher direct images of reflexive Kähler differentials. More precisely, if we have a proper toric morphism $\pi : Y \to X$ from a simplicial toric variety $Y$, the computation of $R^{i}\pi\lsta \Omega_{Y}\ubind{p}$ essentially boils down to a linear algebra computation of finite dimensional vector spaces. Moreover, if $\pi : Y \to X$ is a birational toric morphism such that the polytopal complex associated to $Y$ has every face shellable, then the first assertion of Proposition \ref{prop-formulaKahlerdiff} holds in that setting as well. It would be interesting to investigate what happens for the non-shellable subdivisions.
    \end{rema}
    
    \section{Graded de Rham complex of the Intersection Cohomology Hodge module} \label{section:grDR}

    In this section, we define the generating function associated to the graded de Rham complex of the intersection cohomology Hodge module and prove the main result of the paper.
    
    For $\mu, \tau \in \Sigma_{X}$ such that $\mu \subset \tau$, we define $\dR_{\mu, \tau}$ as
        $$ \dR_{\mu, \tau}(K, L) = \sum_{k, l} \dim_{\CC} \cH^{l} (\gr_{k} \DR_{X} \IC_{S_{\mu}}^{H})_{u} K^{k}L^{l} $$
    for $u \in \tau_{\circ}\sta$. This definition is independent of the choice of $u\in \tau_{\circ}\sta$, as the next result shows.

    \begin{theo} \label{theorem:main-theorem-grDR}
        Let $X$ be the affine toric variety associated to a full dimensional cone $\sigma$ of dimension $n$ and let $\Sigma_{X}$ be the associated fan. Then $\dR_{\mu, \tau}$ is related to $\widetilde{H}_{\mu, \tau}$ in the following way:
        $$ \dR_{\mu, \tau}(K, L) = \widetilde{H}_{\mu, \tau}(K^{-\frac{1}{2}}L) K^{\frac{d_{\mu}-d_{\tau}}{2}} (K^{-1} + L^{-1})^{n - d_{\tau}}. $$
        In particular, $\dR_{\mu, \tau}$ depends only on the graded poset structure of $\Sigma_{X}$. Moreover, $\dR_{\mu,\tau}$ can be explicitly computed in terms of the combinatorics of $\Sigma_X$.
    \end{theo}

    Before giving the proof, we state a lemma relating the graded de Rham complex of the intersection cohomology and the pushforward of Kähler differentials.

    \begin{lemm} \label{lemma:deRham-vs-Kahler}
        Let $\pi \colon Y \to X$ be a birational toric morphism given by a simplicial subdivision of the fan $\Sigma_{X}$. Then for each $\tau \in \Sigma_{X}$, we have
        $$ \Omega_{\tau}(K, L) = \sum_{0 \subset \mu \subset \tau} \dR_{\mu, \tau}(K, L) \cdot D_{\mu}(L^{-1}K^{\frac{1}{2}}) K^{-\frac{d_{\mu}}{2}}.$$
    \end{lemm}
    \begin{proof}
        This is a simple consequence of Proposition \ref{prop:saito-dec-thm-toric}: 
         $$  \pi\lsta \IC_{Y}^{H} = \bigoplus_{\mu \subset \sigma} \bigoplus_{j} \left( \IC_{S_{\mu}}^{H} ( - \frac{d_{\mu} + j}{2} )\right)^{\oplus s_{\mu, j}} [-j].$$
        Note from Remark \ref{rema:Basics-Toric-varieties}. (3) that
         $$ \gr_{k} \DR_Y \IC_{Y}^{H} = \gr_{k} \DR_Y \QQ_{Y}^{H}[n] = \Omega_{Y}\ubind{-k}[n+k]. $$
         By taking $\cH^{l}\gr_{k} \DR$ and using $\bfR f\lsta \circ \gr \DR \simeq \gr \DR \circ  f\lsta$ (Equation \ref{equation-grDRcommuteswithpropermap}), we get
        $$ R^{n+k+l}\pi\lsta (\Omega_{Y}\ubind{-k}) \simeq \cH^{l}\gr_{k} \DR_X \IC_{X}^{H} \oplus \bigoplus_{0 \neq \mu \subset \sigma} \bigoplus_{j} \cH^{l-j} \gr_{k + \frac{d_{\mu} + j}{2}} \DR_{X}(\IC_{S_{\mu}}^{H})^{\oplus s_{\mu, j}}.$$
        By Proposition \ref{prop-formulaKahlerdiff} and induction on dimension, we see that the dimension of the degree $u$ piece of $\cH^{l}\gr_{k} \DR_X \IC_{X}^{H}$ does not depend on the choice of $u \in \tau_{\circ}\sta$. By taking the degree $u$ piece for $u \in \tau_{\circ}\sta$, we get
        $$ \Omega_{\tau}(K, L) = \dR_{0, \tau}(K, L) + \sum_{0 \neq \mu\subset \tau} \dR_{\mu,\tau}(K, L) D_{\mu}(LK^{-\frac{1}{2}}) K^{-\frac{d_{\mu}}{2}}.$$
        The assertion of the lemma follows because we have $D_{\mu}(q) = D_{\mu}(q^{-1})$ by Poincar\'e duality.
    \end{proof}
    
    Now, we give the proof of Theorem \ref{theorem:main-theorem-grDR}.
    
    \begin{proof}[Proof of Theorem \ref{theorem:main-theorem-grDR}]
        We prove this by induction on the dimension of $X$. If $X$ is of dimension zero, then there is nothing to prove. For $\mu \neq 0$, we have the equality
        $$ \dR_{\mu, \tau}^{X} = \dR_{0, \bar{\tau}_{\mu}}^{S_{\mu}} $$
        by definition and the description of the fan of $S_{\mu}$. Note that we also have
        $$ d_{\tau} - d_{\mu} = d_{\bar{\tau}_{\mu}} - d_{0}, \qquad \dim X - d_{\tau} = \dim S_{\mu} - d_{\bar{\tau}_{\mu}}. $$
        Hence, the equality follows by Lemma \ref{lemma-ICstalksforToricSubvar} and the induction hypothesis. Therefore, it is enough to show the equality when $\mu = 0$. Consider the proper toric morphism $\pi \colon Y \to X$ induced by the barycentric subdivision of $\Sigma_{X}$. By Lemma \ref{lemma:deRham-vs-Kahler} and the inductive hypothesis, we have
        \begin{align*}
           \Omega_{\tau}(K, L) & = \dR_{0, \tau}(K, L) + \sum_{0 \neq \mu \subset \tau} \dR_{\mu, \tau}(K, L) \cdot D_{\mu}(L^{-1}K^{\frac{1}{2}}) K^{-\frac{d_{\mu}}{2}}  \\
           & = \dR_{0, \tau}(K, L) + K^{-\frac{d_{\tau}}{2}}(K^{-1} + L^{-1})^{n-d_{\tau}} \sum_{0 \neq \mu \subset \tau} \widetilde{H}_{\mu, \tau}(LK^{-\frac{1}{2}}) D_{\mu}(L^{-1} K^{\frac{1}{2}}).
        \end{align*}
        By Proposition \ref{prop-formulaKahlerdiff}, we have
        \begin{align*}
            \Omega_{\tau}(K, L) & = L^{-n}(1 + K^{-1}L)^{n-d_{\tau}} F_{\tau} (LK^{-\frac{1}{2}}) \\
            & = (K^{-1}+ L^{-1})^{n-d_{\tau}} K^{-\frac{d_{\tau}}{2}} \widetilde{F}_{\tau}(LK^{-\frac{1}{2}}).
        \end{align*}
        Equation \ref{equation:stalk-decomposition-tilde} in Section \ref{section-Decompositiontheorem} gives
        $$ \widetilde{F}_{\tau}(LK^{-\frac{1}{2}}) = \widetilde{H}_{0, \tau}(LK^{-\frac{1}{2}}) + \sum_{0 \neq \mu \subset \tau} \widetilde{H}_{\mu, \tau}(LK^{-\frac{1}{2}}) D_{\mu}(LK^{-\frac{1}{2}}).$$
        By multiplying $K^{-\frac{d_{\tau}}{2}} (K^{-1} + L^{-1})^{n-d_{\tau}}$ on both sides, we get
        $$ \dR_{0, \tau}(K, L) = \widetilde{H}_{0, \tau}(LK^{-\frac{1}{2}}) K^{-\frac{d_{\tau}}{2}} (K^{-1} + L^{-1})^{n-d_{\tau}}.$$
        
        Finally, observe that Remark \ref{rema:explicit-computation-of-H} applied to the barycentric subdivision $\pi$ tells us that $\widetilde{H}_{\mu,\tau}$ can be explicitly computed in terms of the combinatorics of $\Sigma_X$. Therefore, $\dR_{\mu, \tau}$ can be explicitly computed in terms of the combinatorics of $\Sigma_X$ as well.
        \end{proof}

    \begin{rema}\label{rema:relation-to-MS}
        We end the section by relating Theorem \ref{theorem:main-theorem-grDR} to a recent $K$-theoretic result of Maxim and Sch\"urmann. Roughly speaking, the graded de Rham complex gives a homomorphism from $K_{0}^{\mathbb{T}}(\MHM(X))$ to $K_{0}^{\mathbb{T}}(X)[y^{\pm 1}]$, and one can consider the image of $\IC_{X}^{H}[-n]$ by this map. \cite{maxim2024weighted}*{Corollary 5.3} says that the image can be written as the sum
        \[ \sum_{\tau \in \Sigma_X} \chi_{y}(\IC_{X}^{H}[-n] |_{x_{\tau}}) \cdot (1+y)^{n-d_\tau} \cdot (k_\tau)_* [\omega_{S_\tau}]_{\mathbb{T}}.  \]
        Theorem \ref{theorem:main-theorem-grDR} applied to $\dR_{0,\tau}$ gives:
        \begin{align*}
            \dR_{0, \tau}(K, L) &= \widetilde{H}_{0, \tau}(LK^{-\frac{1}{2}}) K^{-\frac{d_{\tau}}{2}} (K^{-1} + L^{-1})^{n-d_{\tau}}\\
            &= \widetilde{H}_{0, \tau}(LK^{-\frac{1}{2}}) (LK^{-\frac{1}{2}})^{d_{\tau}} L^{-d_{\tau}} (K^{-1} + L^{-1})^{n-d_{\tau}}.
        \end{align*}
        Since we are taking the image in $K_{0}^{\bb{T}}(X)[y^{\pm 1}]$, we specialize to $L = -1$ and set $q = LK^{-\frac{1}{2}}$ and $q^2 = K^{-1} = -y$ to get
        \begin{align*}
            \dR_{0, \tau}((-y)^{-1}, -1) &= \widetilde{H}_{0, \tau}(q) q^{d_{\tau}} (1+y)^{n-d_{\tau}}(-1)^n.
        \end{align*}
        We observe that $\chi_{y}(\IC_{X}^{H}[-n]|_{x_{\tau}}) = \widetilde{H}_{0, \tau}(q) q^{d_{\tau}}$ to get
        \begin{align*}
            \dR_{0, \tau}((-y)^{-1}, -1) &= \chi_{y}(\IC_{X}^{H}[-n]|_{x_{\tau}}) (1+y)^{n-d_{\tau}}(-1)^n.
        \end{align*}
        The $(-1)^n$ comes from the fact that $\dR_{0,\tau}$ is defined for $\IC^H_X$ while \cite{maxim2024weighted} work with $\IC^H_X[-n]$.
    \end{rema}

    \appendix
    \section{Explicit Formulas}
    In the appendix, we demonstrate that the polynomials $\widetilde{H}_{\mu, \tau}(q)$ and $\dR_{\mu, \tau}(K, L)$ can be calculated rather explicitly by computing them for full dimensional affine toric varieties up to dimension 4.
    \subsection{Dimension 0, 1, and 2}
    Note that up to dimension 2, every toric variety is simplicial. Hence the intersection cohomology Hodge module agrees with the trivial one. For dimension zero, a zero dimensional toric variety is just a point. Hence 
    $$\widetilde{H}_{0, 0} = q^{0} , \qquad \text{and} \qquad \dR_{0, 0} = K^{0}L^{0}.$$
    For dimension 1, we denote the nonzero ray by $\sigma$. It is easy to check that
    $$ \widetilde{H}_{0, \sigma} = q^{-1}, \quad \widetilde{H}_{0, 0} = \widetilde{H}_{\sigma, \sigma} = q^{0}.$$
    Hence, we get
    $$ \dR_{0, \sigma} = L^{-1}, \quad \dR_{0, 0} = (K^{-1} + L^{-1}), \quad \dR_{\sigma, \sigma} = K^{0}L^{0}.$$

    For dimension 2, let $\tau_{1}, \tau_{2}$ be the two extremal rays of the two dimensional cone $\sigma$. Then, it is clear that
    $$ \widetilde{H}_{0, \sigma} = q^{-2},\quad \widetilde{H}_{0, \tau_{1}} = \widetilde{H}_{0, \tau_{2}} = q^{-1}, \quad  \widetilde{H}_{0, 0} = q^{0}. $$
    The case when the first index is non-zero is redundant since it comes from a lower dimensional toric variety. Therefore,
    $$ \dR_{0, \sigma} = L^{-2} , \quad \dR_{0, \tau_{1}} = \dR_{0, \tau_{2}} = L^{-1}(K^{-1}+L^{-1}), \quad \dR_{0, 0} = (K^{-1} + L^{-1})^{2}.$$

    \subsection{Dimension 3}
    Let $X$ be a 3-dimensional affine toric variety corresponding to a cone $\sigma$ and suppose there are $v$ extremal rays $\mu_{1},\ldots, \mu_{v}$. Note that the number of two dimensional faces is also $v$. Let $\tau_{1},\ldots, \tau_{v}$ be the two dimensional faces. Adding a ray $\rho$ in the interior of the cone gives a proper birational toric morphism $\pi \colon Y \to X$ where $Y$ is simplicial. Note that $\pi$ is an isomorphism outside of the torus fixed point $x_{\sigma}$. Hence, the decomposition theorem tells us
    $$ D_{0}(q) = q^{0}, \quad  D_{\mu_{i}}(q) = D_{\tau_{i}}(q) = 0.$$
    It remains to calculate $D_{\sigma}(q)$. Using (\ref{equation:stalk-decomposition-tilde}) and Proposition \ref{prop-fibercohomology}, we get
    $$ q^{-3}\big((q^{2}-1)^{2} + v(q^{2}-1) + v\big) = \widetilde{H}_{0, \sigma}(q) + D_{\sigma}(q).$$
    Then we can conclude that we have
    $$ D_{\sigma}(q) = q^{1} + q^{-1}, \quad \text{and} \quad \widetilde{H}_{0, \sigma}(q) = q^{-3} + (v-3)q^{-1}.$$
    Note that all the other information for $\widetilde{H}$ comes from lower-dimensional toric varieties because of Lemma \ref{lemma-ICstalksforToricSubvar}.
    Therefore,
    \begin{align*}
        \dR_{0, \sigma} & = L^{-3} + (v-3)K^{-1}L^{-1} \\
        \dR_{0, \tau_{i}} & = (K^{-1} + L^{-1})L^{-2}\\
        \dR_{0, \mu_{i}} & = (K^{-1} + L^{-1})^{2} L^{-1}\\
        \dR_{0, 0} & = (K^{-1} + L^{-1})^{3}.
    \end{align*}

    \subsection{Dimension 4.} Let $X$ be a 4-dimensional affine toric variety corresponding to a cone $\sigma$. We denote the 1-dimensional faces by $\mu_{1},\ldots, \mu_{v}$, the 2-dimensional faces by $\nu_{1},\ldots, \nu_{e}$, and the 3-dimensional faces by $\tau_{1},\ldots, \tau_{f}$. For each $\tau_{k}$, we let $n_{k}$ the number of 1-dimensional faces contained in $\tau_{k}$. We add rays $\rho_{\circ}$ and $\rho_{1}, \ldots, \rho_{f}$ in the interior of $\sigma$ and the interior of $\tau_{k}$'s. This gives a proper birational toric morphism $\pi \colon Y \to X$ where $Y$ is simplicial. The morphism $\pi$ is an isomorphism outside a dimension 1 subset of $X$ and this implies
    $$ D_{0}(q) = q^{0}, \quad D_{\mu_{i}} = D_{\nu_{j}} = 0.$$
    Hence, it remains to calculate $D_{\tau_{k}}$'s and $D_{\sigma}(q)$. By considering the fiber $\pi^{-1}(x_{\tau_{k}})$, we can see that
    $$ D_{\tau_{k}}(q) = q + q^{-1}, \quad \widetilde{H}_{0, \tau_{k}}(q) = q^{-3} + (n_{k} -3)q^{-1}.$$
    This is exactly the same computation as in the previous section. Now, we compute the cohomology of the fiber $\pi^{-1}(x_{\sigma})$. Notice that $d_{1}(\sigma) = 1$ coming from $\rho_{\circ}$. Also,
    $$ d_{2}(\sigma) = \sum_{k} d_{1}(\tau_{k}) + \sum_{j}d_{1}(\nu_{j}) + \sum_{i} d_{1}(\mu_{i}) = f + v = e + 2$$
    since $d_{1}(\tau_{k}) = 1$ coming from $\rho_{i}$ and $d_{1}(\mu_{i}) = 1$ coming from $\mu_{i}$ itself. The last equality follows from Euler's identity $v - e + f = 2$. Also,
    $$ d_{3}(\sigma) = \sum_{k} d_{2}(\tau_{k}) + \sum_{j} d_{2} (\nu_{j}) = \sum_{k} \left( \sum_{\nu_{j} \subset \tau_{k}} d_{1}(\nu_{j}) + \sum_{\mu_{i} \subset \tau_{k}} d_{1}(\mu_{i}) \right) + \sum_{j} d_{2}(\nu_{j}) = \sum_{k} n_{k} + e. $$
    We notice that $n_{k}$ equals the number of 2-dimensional faces of $\tau_{k}$, and for each 2-dimensional face $\mu_{j}$, there are exactly two 3-dimensional faces containing $\mu_{j}$. Therefore
    $$ \sum_{k} n_{k} = 2e.$$
    Similarly,
    $$ d_{4}(\sigma) = \sum_{k} d_{3}(\tau_{k}) = \sum_{k} \sum_{\nu_{j} \subset \tau_{k}} d_{2}(\nu_{j}) = \sum_{k} n_{k} = 2e.$$
    Using (\ref{equation-stalkdecomposition}) and Proposition \ref{prop-fibercohomology}, we get
    \begin{align*}
        & \widetilde{H}_{0, \sigma}(q) + D_{\sigma}(q)\\
        & = q^{-4}((q^{2}-1)^{3} + (e + 2) (q^{2}-1)^{2} + 3e(q^{2}-1) + 2e) - \sum_{k} \widetilde{H}_{\tau_{k}, \sigma}(q) D_{\tau_{k}}(q) \\
        & = q^{-4}(q^{6} + (e-1)q^{4} + (e-1)q^{2} + 1) - f \cdot (q + q^{-1})q^{-1} \\
        & = q^{2} + (e-1-f)q^{0} + (e-1-f)q^{-2} + q^{-4} \\
        & = q^{2} + (v-3)q^{0} + (v-3)q^{-2} + q^{-4}.
    \end{align*}
    Therefore, $D_{\sigma}(q) = q^{2} + (v-3)q^{0} + q^{-2}$ and
    $$ \widetilde{H}_{0, \sigma}(q) = q^{-4} + (v-4)q^{-2}.$$
    Therefore,
    \begin{align*}
        \dR_{0, \sigma} & = L^{-4} + (v-4) K^{-1} L^{-2} \\
        \dR_{0, \tau_{k}} & = (K^{-1} + L^{-1}) \cdot (L^{-3} + (n_{k}-3) K^{-1}L^{-1}) \\
        \dR_{0, \nu_{j}} & = (K^{-1} + L^{-1})^{2} \cdot L^{-2} \\
        \dR_{0, \mu_{i}} & = (K^{-1} + L^{-1})^{3} \cdot L^{-1} \\
        \dR_{0, 0} & = (K^{-1} + L^{-1})^{4}.
    \end{align*}

    \section{A reinterpretation of Theorem \ref{theo:main-theo1.1}} \label{appendix-2}
    Let $X$ be an affine toric variety corresponding to a strictly convex polyhedral full dimensional cone $\sigma$ of dimension $n$. We first prove a couple of lemmas:

    \begin{lemm} \label{lemm:strata-of-IC-pure-HT}
        Let $\tau$ be a face of $\sigma$ and let $j_{\tau} \colon O_{\tau} \to X$ be the inclusion. Then $j_{\tau}\sta \IC_{X}^{H}$ is a pure complex of weight $n$, and
        $$ j_{\tau}\sta \IC_{X}^{H} \simeq \bigoplus_{l} \cH^{l}(\IC_{X})_{x_{\tau}} \otimes \QQ_{O_{\tau}}^{H} \left( -\frac{n + l}{2}\right) [-l].$$
	\end{lemm}
	\begin{proof} It is enough to show that $\cH^{j} (i_{x_{\sigma}}\sta \IC_{X}^{H})$ is pure of Hodge--Tate type with weight $j +n$, where $x_{\sigma}$ is the torus fixed point. Then we can use the strategy in the proof of Proposition \ref{prop:saito-dec-thm-toric} to get the corresponding statement about the other orbits, by writing the affine open $U_{\tau}$ as $ U_{\tau} \simeq O_{\tau} \times V_{\tau}$. 

		Consider the morphism $\pi \colon \widetilde{X} \to X$ given by a barycentric subdivision of $\sigma$ as in Remark \ref{rema:barycentric-resolution}. Let $E = \pi^{-1}(x_{\sigma})$ be the inverse image of the torus fixed point. We point out that $E$ is an irreducible proper simplicial toric variety of dimension $n-1$. Note that $\IC_{X}^{H}$ is a summand of $\pi\lsta \QQ_{\widetilde{X}}^{H}[n]$ by the Decomposition theorem. By proper base change (see \cite{saito1990mixedHodgemodules}*{4.4.3}), we see that $i_{x_{\sigma}}\sta \IC_{X}^{H}$ is a direct summand of $\pi\lsta \QQ_{E}^{H}[n]$. Therefore, $\cH^{l}i_{x_{\sigma}}\sta \IC_{X}^{H}$ is a direct summand of $H^{n+l}(E, \QQ_{E})$ which is pure of Hodge--Tate type of weight $n+l$, by \cite{dCMM-toricmaps}*{Theorem A}.

        The purity of the complex tells us that $j\sta_{\tau} \IC^H_{X}$ decomposes in the derived category of mixed Hodge modules.
	\end{proof}

    \begin{lemm} \label{lemm:j!-of-orbits}
        Let $j_{\tau} \colon O_{\tau} \to X$ as in Lemma \ref{lemm:strata-of-IC-pure-HT} and let $S_{\tau}$ be the torus invariant closed subvariety corresponding to $\tau$. Then
        $$ \gr_{-k} \DR_{X} j_{!} \QQ_{O_{\tau}}^{H} \simeq \omega_{S_{\tau}}[-k]^{\oplus {\dim O_{\tau} \choose k}}.$$
    \end{lemm}
    \begin{proof}
        Consider a toric resolution $\pi \colon \widetilde{S}_{\tau} \to S_{\tau}$ and denote by $\widetilde{j} \colon O_{\tau} \to \widetilde{S}_{\tau}$ the open inclusion and $\widetilde{D}_{\tau}$ its complement. First, we show that
        $$ \gr_{-k} \DR_X j\lsta \QQ_{O_{\tau}}^{H} \simeq \cO_{S_{\tau}}^{\oplus {\dim O_{\tau} \choose k} }[-k].$$
        This follows from
        $$ \gr_{-k} \DR_X \widetilde{j}\lsta \QQ_{O_{\tau}}^{H} \simeq \Omega_{\widetilde{S}_{\tau}}^{k}(\log D_{\tau}) [-k] \simeq O_{\widetilde{S}_{\tau}}^{\oplus {\dim O_{\tau} \choose k}}[-k], $$
        the commutativity of the graded de Rham complex with proper push-forward, and the fact that toric varieties have rational singularities. From $j_{!} \QQ_{O_{\tau}}^{H} \simeq \DD j\lsta \DD \QQ_{O_{\tau}}^{H}$ and the commutativity between duality of Hodge modules and Grothendieck duality, we get
        \begin{align*}
            &\gr_{-k} \DR \DD j_{\ast} \DD \QQ_{O_{\tau}}^{H} \\
            & \simeq \bfR \SheafHom_{\cO_{S_{\tau}}}(\gr_{k} \DR j\lsta \QQ_{O_{\tau}}^{H}[2\dim O_{\tau}](\dim O_{\tau}), \omega_{S_{\tau}}[\dim O_{\tau}])\\
            & \simeq \bfR\SheafHom_{\cO_{S_{\tau}}}(\gr_{-\dim O_{\tau} + k} \DR j\lsta \QQ_{O_{\tau}}^{H}, \omega_{S_{\tau}} ) [-\dim O_{\tau}] \\
            & \simeq \bfR\SheafHom_{\cO_{S_{\tau}}}(\cO_{S_{\tau}}[-\dim O_{\tau}+k] , \omega_{S_{\tau}})^{\oplus {\dim O_{\tau}\choose \dim O_{\tau} - k}} [-\dim O_{\tau}] \simeq  \omega_{S_{\tau}}^{\oplus {\dim O_{\tau} \choose k}}[-k].
        \end{align*}
    \end{proof}

    We now give an alternate proof of Theorem \ref{theo:main-theo1.1}.

    \begin{proof} We point out that as in the beginning of the proof of Theorem \ref{theorem:main-theorem-grDR}, it is enough to verify the equation concerning $\dR_{0, \tau}$, where $\tau$ is a face of $\sigma$.
        
        Let $X_{m}$ be the union of all $m$-dimensional torus invariant closed subvarieties of $X$ and $i_{m} \colon X_{m} \to X$ be the closed immersion. Note that $X_{m} \setminus X_{m-1}$ is a disjoint union of all $m$-dimensional orbits. We denote by $j_{m} \colon X_{m} \setminus X_{m-1} \to X$ the locally closed embedding. From \cite{saito1990mixedHodgemodules}*{4.4.1}, we have an exact triangle
    $$ \to j_{m!} j_{m}\sta \IC_{X}^{H} \to i_{m}\sta \IC_{X}^{H} \to i_{m-1}\sta \IC_{X}^{H} \xrightarrow{+1}.$$
        We show by induction on $m$ that 
        \begin{enumerate}
            \item $\cH^{l} \gr_{k} \DR i_{m}\sta \IC^H_{X}$ is supported in degrees $u \in \tau_{\circ}\sta$ such that $\dim \tau \geq n-m$, and
            \item the induced long exact sequence of graded de Rham complexes
            $$ \cdots \to \cH^{l} \gr_{k} \DR j_{m!}j_{m}\sta \IC_{X}^{H} \to \cH^{l} \gr_{k} \DR i_{m}\sta \IC_{X}^{H} \to \cH^{l} \gr_{k} \DR i_{m-1}\sta \IC_{X}^{H} \to \cdots $$
            splits into short exact sequences
             $$ 0 \to \cH^{l} \gr_{k}\DR j_{m!}j_{m}\sta \IC_{X}^{H} \to \cH^{l} \gr_{k}\DR i_{m}\sta \IC_{X}^{H} \to \cH^{l} \gr_{k} \DR i_{m-1}\sta \IC_{X}^{H} \to 0. $$
        \end{enumerate}
        These assertions are clear when $m = 0$. By the induction hypothesis, $\cH^{l} \gr_{k}\DR i_{m-1}\sta \IC_{X}^{H}$ is supported in degrees $u \in \tau_{\circ}\sta$ such that $\dim \tau \geq n-m+1$. Using Lemma \ref{lemm:strata-of-IC-pure-HT} and \ref{lemm:j!-of-orbits} and the properties of $\omega_{S_\tau}$, we conclude that $\cH^{l} \gr_{k} \DR j_{m!}j_{m}\sta \IC_{X}^{H}$ is supported in degrees $u \in \tau_{\circ}\sta$ such that $\dim \tau = n-m$. Therefore, the torus-equivariant homomorphism $\cH^{l} \gr_{k} \DR i_{m-1}\sta \IC_{X}^{H} \to \cH^{l+1} \gr_{k} \DR j_{m!} j_{m}\sta \IC_{X}^{H}$ is zero. This shows both the assertions and so, we see that $\cH^{l} \gr_{k} \DR \IC_{X}^{H}$ admits a filtration whose successive quotients are $\cH^{l} \gr_{k} \DR j_{m!} j_{m}\sta \IC_{X}^{H}$. Then again by Lemma \ref{lemm:strata-of-IC-pure-HT} and \ref{lemm:j!-of-orbits}, we see that for $u \in \tau_{\circ}\sta$, we get
        \begin{align*}
            & \dR_{0, \tau}(K, L) =  \sum_{k, l} \dim \cH^{l} \left( \gr_{k} \DR \IC_{X}^{H}\right)_{u} K^{k} L^{l} \\
            & = \sum_{k, l} \dim \cH^{l} (\IC_{X})_{x_{\tau}} {n - d_{\tau} \choose k - \frac{n+l}{2}} K^{-k} L^{k + \frac{l-n}{2}} \\
            & = \sum_{k,l}\dim \cH^{l} (\IC_{X})_{x_{\tau}} {n - d_{\tau} \choose k - \frac{n+l}{2}} (K^{-1}L)^{k - \frac{l}{2}} (K^{-\frac{1}{2}}L)^{l} L^{-\frac{n}{2}} \\
            & = (1 + K^{-1}L)^{n-d_{\tau}}K^{-\frac{n}{2}} H_{0,\tau} (K^{-\frac{1}{2}}L) \\
            & = (1 + K^{-1}L)^{n - d_{\tau}} K^{-\frac{n}{2}} (K^{-\frac{1}{2}}L)^{-n + d_{\tau}} \widetilde{H}_{0, \tau}(K^{-\frac{1}{2}}L) \\
            & = (K^{-1} + L^{-1})^{n-d_{\tau}} K^{-\frac{d_{\tau}}{2}} \widetilde{H}_{0, \tau}(K^{-\frac{1}{2}}L)
        \end{align*}
    \end{proof}

    {\bf Acknowledgments.} We would like to thank Mircea Musta\c{t}\u{a} for numerous helpful discussions, and Claudiu Raicu for pointing out the reference \cite{CFS-Effectivedecompositiontheorem} which helped us compute the polynomial $\widetilde H_{\mu,\tau}$ explicitly, and J\"org Sch\"urmann for kindly explaining the relation between \cite{maxim2024weighted} and Theorem \ref{theo:main-theo1.1}. We would also like to thank the referee for a very helpful suggestion of an alternate approach to Theorem \ref{theo:main-theo1.1}.
    
    \bibliographystyle{alpha}
    \bibliography{Reference}

\end{document}